\documentclass{article}

\usepackage[nonatbib, preprint]{neurips_2024}

\usepackage[utf8]{inputenc} % allow utf-8 input
\usepackage[T1]{fontenc}    % use 8-bit T1 fonts
\usepackage{hyperref}       % hyperlinks
\usepackage{url}            % simple URL typesetting
\usepackage{booktabs}       % professional-quality tables
\usepackage{amsfonts}       % blackboard math symbols
\usepackage{amssymb}
\usepackage{amsthm}
\usepackage{nicefrac}       % compact symbols for 1/2, etc.
\usepackage{microtype}      % microtypography
\usepackage{xcolor}         % colors
\usepackage{graphicx}
\usepackage{float}
\usepackage{multirow}
\usepackage{tabularx}
\usepackage{subcaption}
\usepackage{enumitem}

\usepackage[round]{natbib}
\usepackage{amsmath,amsfonts,bm}

\def\eqref#1{equation~\ref{#1}}
\def\floor#1{\lfloor #1 \rfloor}
\def\1{\bm{1}}

\def\vb{{\bm{b}}}
\def\vc{{\bm{c}}}

\def\ve{{\bm{e}}}

\def\vh{{\bm{h}}}

\def\vl{{\bm{l}}}

\def\vp{{\bm{p}}}

\def\vr{{\bm{r}}}

\def\vu{{\bm{u}}}
\def\vv{{\bm{v}}}
\def\vw{{\bm{w}}}
\def\vx{{\bm{x}}}
\def\vy{{\bm{y}}}
\def\vz{{\bm{z}}}

\def\mA{{\bm{A}}}

\def\mD{{\bm{D}}}

\def\mH{{\bm{H}}}
\def\mI{{\bm{I}}}

\def\mN{{\bm{N}}}

\def\mP{{\bm{P}}}
\def\mQ{{\bm{Q}}}

\def\mT{{\bm{T}}}

\def\mW{{\bm{W}}}
\def\mX{{\bm{X}}}
\def\mY{{\bm{Y}}}
\def\mZ{{\bm{Z}}}

\DeclareMathAlphabet{\mathsfit}{\encodingdefault}{\sfdefault}{m}{sl}
\SetMathAlphabet{\mathsfit}{bold}{\encodingdefault}{\sfdefault}{bx}{n}

\newcommand{\R}{\mathbb{R}}

\newtheoremstyle{theoremstyle} % name
    {1em}                % Space above
    {1em}                % Space below
    {}                   % Body font
    {}                   % Indent amount
    {\bfseries}          % Theorem head font
    {.}                  % Punctuation after theorem head
    {.5em}               % Space after theorem head
    {\thmname{#1}\thmnumber{ #2}\thmnote{ ({{\bfseries #3}})}} % Theorem head spec

\theoremstyle{theoremstyle}
\newtheorem{theorem}{Theorem}
\newtheorem{lemma}[theorem]{Lemma}
\newtheorem{definition}{Definition}

\usepackage{algorithm}
\usepackage{algorithmic}

\title{NeuralQP: A General Hypergraph-based Optimization Framework for Large-scale QCQPs}

\author{
  Zhixiao Xiong\thanks{Department of Computer Science and Technology, Tsinghua University}\;\,\thanks{Weiyang College, Tsinghua University} \\
  \texttt{xiong-zx21@mails.tsinghua.edu.cn}
  \And
  Fangyu Zong\footnotemark[1]\;\,\thanks{Department of Industrial Engineering, Tsinghua University}\\
  \texttt{zongfy21@mails.tsinghua.edu.cn}
  \And
  Huigen Ye\footnotemark[1]\\
  \texttt{yhg23@mails.tsinghua.edu.cn}
  \And
  Hua Xu\footnotemark[1]\;\,\thanks{Corresponding Author}\\
  \texttt{xuhua@tsinghua.edu.cn}
}

\begin{document}

\maketitle

\begin{abstract}
    Machine Learning (ML) optimization frameworks have gained attention for their ability to accelerate the optimization of large-scale Quadratically Constrained Quadratic Programs (QCQPs) by learning shared problem structures. However, existing ML frameworks often rely heavily on strong problem assumptions and large-scale solvers. This paper introduces NeuralQP, a general hypergraph-based framework for large-scale QCQPs. NeuralQP features two main components: Hypergraph-based Neural Prediction, which generates embeddings and predicted solutions for QCQPs without problem assumptions, and Parallel Neighborhood Optimization, which employs a McCormick relaxation-based repair strategy to identify and correct illegal variables, iteratively improving the solution with a small-scale solver. We further prove that our framework UniEGNN with our hypergraph representation is equivalent to the Interior-Point Method (IPM) for quadratic programming. Experiments on two benchmark problems and large-scale real-world instances from QPLIB demonstrate that NeuralQP outperforms state-of-the-art solvers (e.g., Gurobi and SCIP) in both solution quality and time efficiency, further validating the efficiency of ML optimization frameworks for QCQPs.
\end{abstract}

\section{Introduction}
\label{sec:introduction}
Quadratically Constrained Quadratic Programs (QCQPs) are mathematical optimization problems characterized by the presence of quadratic terms, finding extensive applications across diverse domains such as finance \citep{gondzio2007parallel}, robotic control \citep{galloway2015torque}, and power grid operations \citep{zhang2013sequential}. However, solving QCQPs is exceptionally challenging due to their discrete \citep{balas1969duality} and nonconvex nature \citep{elloumi2019global}, especially for large-scale QCQPs. With the advancement of machine learning (ML), the ML-based QCQP optimization framework has emerged as a promising research direction as it can effectively leverage shared problem structures among similar QCQPs to accelerate the solution process.

Current widely adopted ML-based QCQP optimization frameworks can be categorized into two types: \textit{Solver-based Learning} and \textit{Model-based Learning}. \textit{Solver-based Learning} methods learn to tune the parameters or status of the solver to accelerate the solution process.  As a representative work, RLQP \cite{ichnowski2021accelerating} learned a policy to tune parameters of OSQP solver \citep{stellatoOSQPOperatorSplitting2017} with reinforcement learning (RL) to accelerate convergence. \cite{bonamiLearningClassificationMixedInteger2018} learned a classifier that predicts a suitable solution strategy on whether or not to linearize the problem for the CPLEX solver. \cite{ghaddarLearningSpatialBranching2022} and \cite{kannanLearningAccelerateGlobal2023} both learned branching rules on selected problem features to guide the solver. Although solver-based methods have demonstrated strong performance on numerous real-world QCQPs, their effectiveness heavily relies on large-scale solvers and is constrained by the solver's solving capacity, resulting in scalability challenges.

\textit{Model-based learning} methods employ a neural network model to learn the parameters of the QCQP models, aiming to translate the optimization problem into a multiclass classification problem, of which the results can accelerate the solution process of the original optimization problem.
\cite{bertsimasVoiceOptimization2020, bertsimasOnlineMixedIntegerOptimization2021} learned a multi-class classifier for both solution strategies and integer variable values, proposing an online QCQP optimization framework that consists of a feedforward neural network evaluation and a linear system solution. However, these methods make strong assumptions about the model parameters. They assume that multiple problem instances are generated by a parametric model and that the parameters can be correctly inferred from data, which limits their generalization ability to more complicated problems.

To address the above limitations, this paper proposes NeuralQP, a general hypergraph-based optimization framework for large-scale QCQPs, which can be divided into two stages. 
In \textit{Hypergraph-based Neural Prediction}, a representation of QCQPs based on a variable relational hypergraph with initial vertex embeddings is created as a complete representation of QCQPs.Then UniEGNN takes both hyperedge features and vertex features as inputs, leverages an enhanced vertex-hyperedge-vertex convolution strategy, and finally obtains neural embeddings for the variables in QCQPs. In \textit{Parallel Neighborhood Optimization}, to recover feasible solutions after multiple neighborhood solutions are merged (termed crossover), a new repair strategy based on the McCormick relaxation is proposed to quickly identify violated constraints. Then, improperly fixed variables are reintroduced into the neighborhood, which realizes an adaptive updating of the neighborhood radius and effective correction of infeasible solutions via a small-scale optimizer. 

To validate the effectiveness of NeuralQP, experiments are conducted on two benchmark QCQPs and large-scale instances from QPLIB \citep{furiniQPLIBLibraryQuadratic2019}. Results show that NeuralQP can achieve better results than the state-of-the-art solver Gurobi and SCIP in a fixed wall-clock time, only with a small-scale solver with 30\% of the original problem size. Further experiments validate that NeuralQP can achieve the same solution quality in less than 10\% of the solving time of Gurobi and SCIP in large-scale QCQPs, which validates the efficiency of our framework in solving QCQPs. Our contributions are concluded as follows:

\begin{enumerate}
    \item We propose NeuralQP, the first general optimization framework for large-scale QCQPs using only small-scale solvers without any problem assumption, shedding light on solving general nonlinear programming using ML-based frameworks.
    
    \item In \textit{Hypergraph-based Neural Prediction}, a new hypergraph-based representation is proposed as a complete representation of QCQPs, and an enhanced hypergraph convolution framework UniEGNN is employed to fully utilize hyperedge features. We theoretically prove that UniEGNN is equivalent to the Interior-Point Method (IPM) with this representation.
    
    \item In \textit{Parallel Neighborhood Optimization}, a new repair strategy based on the McCormick relaxation is proposed for neighborhood search and crossover with small-scale solvers, achieving parallel optimization of neighborhoods.
    
    \item Experiments show that NeuralQP can accelerate the convergence speed with a small-scale optimizer and generalize well to QCQPs with much larger sizes, indicating its potential in addressing large-scale QCQPs. We also conduct ablation experiments to verify the effectiveness of both the Neural Prediction and neighborhood search separately.
\end{enumerate}

\section{Preliminaries}

\subsection{Quadratically-Constrained Quadratic Program}
A Quadratically-Constrained Quadratic Program (QCQP) \citep{elloumi2019global} is an optimization problem that involves minimizing (or maximizing) a quadratic objective function subject to quadratic constraints. Formally, a QCQP instance $I$ is defined as follows:
\begin{equation}
    \begin{aligned}
        \min \text{ or } \max f(\vx) &= \vx^{\mathrm{T}} \mQ^0 \vx+\left(\vr^0\right)^{\mathrm{T}} \vx, \\
        \text {s.t.}\quad & \vx^{\mathrm{T}} \mQ^j \vx+\left(\vr^j\right)^{\mathrm{T}} \vx \le b_j, \quad \forall j \in \mathcal{M},\\
        &l_i \leq x_i \leq u_i, \quad \forall i \in \mathcal{N},\\
        & x_i \in \mathbb{Z}, \quad \forall i \in \mathcal{I}.
    \end{aligned}
    \label{eq:qp}
\end{equation}
where $\mathcal{N}=\{1,2,\dots,n\}$ and $\mathcal{M} = \{1,2,\dots,m\}$ are index sets of variables and constraints, and $\mathcal{I} \subset \mathcal{N}$ is the index set of integer variables. $\vx = (x_1,x_2,\dots,x_n) \in \R ^ n $ denotes the vector of $n$ variables with $l_i$ and $u_i$ being the lower and upper bounds of $x_i$ respectively. For each $j \in \mathcal{M}$, $\mQ^j \in \R ^{n\times n}$ is symmetric but not necessarily positive definite, $\vr^j$ is the coefficient vector of linear terms. For $j \in\mathcal{M} \setminus \{0\}$, $b_j$ denotes the right-hand-side of the $j$-th constraint.

\subsection{Graph Representations for MILPs}
\label{sec:graph_representation}
Graph representations for (mixed-integer) linear programs (denoted as (MI)LPs) have been proposed to transform the MILP into a suitable input format for graph neural networks, accompanied by encoding strategies that further enhance the expressive power of the corresponding graph. This section introduces the specifics of \textit{tripartite graph representation} and \textit{random feature strategy}.

\begin{figure}[t]
    \centering
    \begin{minipage}{0.55\linewidth}
        \includegraphics[width=\linewidth]{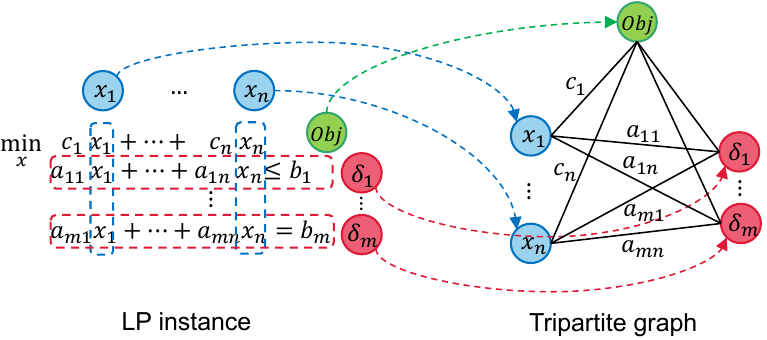}
        \caption{Tripartite representation. 
        The green, blue, and red nodes represent the objective, variables, and constraints, and the edge features are associated with coefficients in the original problem.
        }
        \label{fig:tripartite_graph}
    \end{minipage}
    \hfill
    \begin{minipage}{0.4\linewidth}
        \centering
        \includegraphics[width=\linewidth]{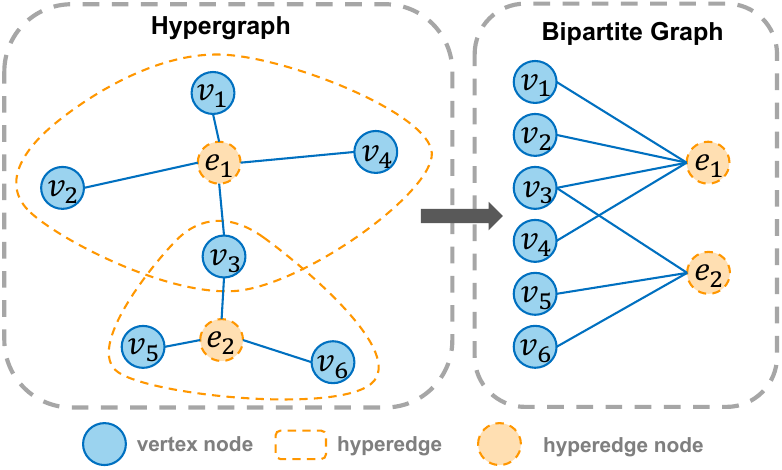}
        \caption{Star expansion. 
        A hypergraph is transformed into a bipartite graph, with nodes on the left and hyperedges on the right.
        }
        \label{fig:hypergraph2bipartite}
    \end{minipage}
\end{figure}

\citet{gasseExactCombinatorialOptimization2019} first proposed a \textit{bipartite graph representation} for MILPs that preserves the entire information of constraints, decision variables, and their relationships without loss. \citet{dingAcceleratingPrimalSolution2019} further proposed a \textit{tripartite graph representation} for MILPs, which simplifies the feature representation of the graph by adding the objective node, ensuring that all coefficients appear only on edge weights. Figure \ref{fig:tripartite_graph} gives an example of the tripartite graph representation. In Figure \ref{fig:tripartite_graph}, the left set of $n$ variable nodes represent decision variables with variable types and bounds encoded into node features; the $m$ constraint nodes on the right symbolize linear constraints with constraint senses and $b_i$ values encoded into node features; the objective node on the top represent the objective function with objective sense encoded as node features. The connecting edge $(i, j)$, signifies the presence of the $i$-th decision variable in the $j$-th constraint, with edge weight $a_{ij}$ representing the coefficient. The upper node signifies the objective, and the edge weight $c_i$ between the $i$-th variable node and the objective represents the objective coefficient. \citet{qian2024exploring} further proved that Message-Passing Neural Networks (MPNNs) with the tripartite graph representation can simulate standard Interior-Point Methods for linear programming.

However, there are MILP instances, named \textit{foldable} MILPs \citep{chenRepresentingMixedIntegerLinear2023} which have distinct optimal solutions yet the corresponding graph cannot be distinguished by the Weisfeiler-Lehman test \citep{weisfeiler1968reduction}. On these foldable MILPs, the performance of the graph neural networks could significantly deteriorate. To resolve such inability, they proposed a \textit{random feature strategy}, i.e., appending an extra dimension of random number to the node features, which further enhances the power of graph neural networks.

\subsection{McCormick Relaxation}
\label{sec:McCormick_relaxation}

The McCormick relaxation \citep{mccormickComputabilityGlobalSolutions1976} is widely employed in nonlinear programming, aiming at bounding nonconvex terms by linear counterparts. Given variables $x$ and $y$ with bounds $L_x \leq x \leq U_x$ and $L_y \leq y \leq U_y$, the nonconvex term $\phi_{xy}:=xy$ can be approximated by:
\begin{equation}
    \begin{aligned}
        \phi_{xy} &\leq \min\{ L_y x+U_x y-L_y U_x,\ L_x y+U_y x-L_x U_y\},\\
        \phi_{xy} &\geq \max \{L_y x+L_x y-L_x L_y,\ U_y x+U_x y-U_x U_y\},
        \label{eq:McCormick}
    \end{aligned}
\end{equation}
which transforms the original nonconvex problem into a more manageable linear format. This method is particularly pertinent when dealing with nonconvex quadratic terms, for which obtaining feasible solutions is notably challenging due to the inherent nonlinearity and nonconvexity.

\subsection{Hypergraph Basics}
\label{sec:hypergraph_basics}

A hypergraph $\mathcal{H} = (\mathcal{V}, \mathcal{E})$ is defined by a set of vertices $\mathcal{V} = \{1,2,\cdots,n\}$ and a set of hyperedges $\mathcal{E} = \left\{e_j\right\}_1^m$, where each hyperedge $e\in \mathcal{E}$ is a non-empty subset of $\mathcal{V}$ \citep{bretto2013hypergraph}. The incidence matrix $\mH \in \{0,1\}^{|\mathcal{V}|\times|\mathcal{E}|}$ is characterized by $\mH(v, e) = 1$ if $v\in e$ else $0$. A hypergraph is termed \emph{$k$-uniform} \citep{rodl2004regularity} if every hyperedge contains exactly $k$ vertices. A distinctive class within hypergraph theory is the \emph{bipartite} hypergraph \citep{annamalai2016finding}. If $\mathcal{V}$ can be partitioned into two disjoint sets $\mathcal{V}_1$ and $\mathcal{V}_2$, a bipartite hypergraph is defined as $\mathcal{H} = (\mathcal{V}_1, \mathcal{V}_2, \mathcal{E})$ so that for each hyperedge $e \in E$, $|e \cap \mathcal{V}_1| = 1$ and $e \cap \mathcal{V}_2 \neq \emptyset$.

Hypergraphs can be converted into graphs via expansion techniques \citep{daiHypergraphComputation2023}. We introduce the \emph{star expansion} technique, which converts a hypergraph into a bipartite graph. Specifically, given a hypergraph $\mathcal{H} = (\mathcal{V}, \mathcal{E})$, each hyperedge $ e \in \mathcal{E} $ is transformed into a new node in a bipartite graph. The original vertices from $ \mathcal{V} $ are retained as nodes on one side of the bipartite graph, while each hyperedge is represented as a node on the opposite side, connected to its neighboring vertices in $ \mathcal{V} $ by edges, thus forming a star. Such transformation enables the representation of the high-order relationships in the hypergraph within the simpler structure of a bipartite graph. Figure \ref{fig:hypergraph2bipartite} provides a bipartite graph representation obtained through the star expansion of a hypergraph.

\begin{definition} \textbf{Inter-Neighbor Relation.}
    The Inter-Neighbor Relation $\mathrm{N} \subset \mathcal{V}\times \mathcal{E}$ on a hypergraph $\mathcal{H} = (\mathcal{V}, \mathcal{E})$ with incidence matrix $\mH$ is defined as: $\mathrm{N} = \{ (v,e) | \mH(v,e) = 1, v \in \mathcal{V}, e \in \mathcal{E}\}.$
    The hyperedge neighborhood $\mathcal{N}_e(v)$ of vertex $v$ and the vertex neighborhood $\mathcal{N}_v(e)$ of hyperedge $e$ are defined based on the Inter-Neighbor Relation.
\end{definition}

\begin{definition}\textbf{Hyperedge Neighborhood.}
    The hyperedge neighborhood of vertex $v \in \mathcal{V}$ is defined as: $\mathcal{N}_e(v) = \{ e | v\mathrm{N} e, e \in \mathcal{E} \}, \text{ for each } v \in \mathcal{V}.$
\end{definition}

\begin{definition}\textbf{Vertex Neighborhood.}
    The vertex neighborhood of hyperedge $e\in \mathcal{V}$ is defined as: $\mathcal{N}_v(e) = \{ v | v\mathrm{N} e, v \in \mathcal{V} \}, \text{ for each }e \in \mathcal{E}.$
\end{definition}

Using hypergraph Inter-Neighbor Relation, the general hypergraph convolution (also known as \textit{spatial convolution}) \citep{gaoHGNNGeneralHypergraph2023} follows a vertex-hyperedge-vertex message propagation pattern. As the first representative framework of hypergraph convolution, \emph{UniGNN} \citep{huangUniGNNUnifiedFramework2021} utilizes a convolution method delineated in Equation \ref{eq:UniGNN}:
\begin{equation}
    \begin{aligned}
        \text { (UniGNN) }\left\{\begin{array}{l}
        h_e=\phi_1\left(\left\{h_j\right\}_{j \in \mathcal{N}_v(e)}\right)\\
        \tilde{h}_v=\phi_2\left(h_v,\left\{h_i\right\}_{i \in \mathcal{N}_e(v)}\right)
        \end{array}\right. ,
        \label{eq:UniGNN}
    \end{aligned}
\end{equation}
where $h_e$ symbolizes the hyperedge features aggregated from the vertices in hyperedge $e$; $h_v$ and $\tilde h_v$ represent the vertex features before and after the convolution respectively; the functions $\phi_1$ and $\phi_2$ are permutation-invariant functions. This spatial convolution strategy initiates by aggregating messages from the incident vertices of a hyperedge and then sends the aggregated message back to vertices, thus fulfilling a round of message propagation.

\section{Hypergraph-based Neural Prediction}
\label{sec:hypergraph_based_neural_prediction}

In this section, we explore the application of hypergraphs to enhance the representation and solution of QCQPs. We first introduce the construction of a \textit{variable relational hypergraph} for QCQP; we then introduce UniEGNN as a hypergraph-based convolution framework for QCQPs and prove its equivalence to the IPM for quadratic programming.

\subsection{Hypergraph-based Representation for QCQP}
\label{sec:hypergraph_based_representation_for_qcqp}

Based on the definition of hypergraph, we propose the \emph{variable relational hypergraph} as a complete representation of QCQPs compatible with MILPs. Consider a QCQP problem defined in Equation \ref{eq:qp}, formal definitions are given below to describe the construction of such a hypergraph.

\begin{definition}\textbf{Extended variable vertex set.} 
    The variable vertex set $\mathcal{V}_x$ is defined as $\mathcal{V}_x = \{v_i|i\in \mathcal{N} \}$ and the extended variable vertex set $\bar{\mathcal{V}}_x$ is defined as $\bar{\mathcal{V}}_x=\mathcal{V}_x \cup \{v_0,v^2\}$, which is the variable vertex set with additional two vertices representing degree zero and degree two.
\end{definition}

\begin{definition}\textbf{Constraint vertex set and objective vertex set.} 
    The constraint vertex set $\mathcal{V}_c$ is defined as $\mathcal{V}_c = \{c_i|i\in \mathcal{M} \}$ and the objective vertex set $\mathcal{V}_o$ is defined as $\mathcal{V}_o = \{o\}$.
\end{definition}

\begin{definition}\textbf{V-O relational hyperedge.}
    The set of V-O relational hyperedges is defined as $\mathcal{E}_o = \{\{v_i, v_0,o\} | \text{ if term }x_i \text{ is in the objective}\} \cup \{\{v_i, v^2,o\} | \text{ if term }x_i^2 \text{ is in the objective}\} \cup \{\{v_i, v_j,o\} | \text{ if term }x_i x_j \text{ is in the objective}\}$.
\end{definition}
\begin{definition}
    \textbf{V-C relational hyperedge.} The set of V-C relational hyperedges is defined as $\mathcal{E}_c = \{\{v_i, v_0,c_k\} | \text{ if term }x_i \text{ is in constraint } k\} \cup \{\{v_i, v^2,c_k\}|\text{ if term }x_i^2 \text{ is in constraint } \delta_k\} \cup \{\{v_i, v_j,c_k\}|\text{ if term }x_i x_j \text{ is in constraint } \delta_k\}$.
\end{definition}
\begin{definition}
    \textbf{Variable relational hypergraph.} A Variable Relational Hypergraph is defined as $\mathcal{H} = (\mathcal{V}_x,\mathcal{V}_o,\mathcal{V}_c,\mathcal{E}_o,\mathcal{E}_c)$ where $\mathcal{V}_x,\mathcal{V}_o,\mathcal{V}_c$ are sets of nodes representing variables, objective, and constraints respectively. $\mathcal{E}_o$ and $\mathcal{E}_c$ are sets of V-O relational and V-C relational hyperedges.
\end{definition}

Each QCQP instance $I$ is transformed into a variable relational hypergraph$\mathcal{H}(I):=(\mathcal{V}_x,\mathcal{V}_o,\mathcal{V}_c,\mathcal{E}_o,\mathcal{E}_c)$. Coefficients of the terms are then encoded as features of the corresponding hyperedges; variable, constraint and objective vertex features are generated in the same way as is introduced in Section \ref{sec:graph_representation}.

As a concrete example shown in Figure \ref{fig:hypergraph_representation}, the QCQP instance is first converted to a variable relational hypergraph with initial embeddings, where the $q_{11}^0 x_1^2$ term in the objective is represented as a hyperedge covering vertices $v_1$, $v^2$ and $o$ with feature $q_{11}^0$; the term $q_{1n}^{1}x_1 x_n$ in constraint $\delta_1$ is represented as a hyperedge covering vertices $v_1$, $v_0$ and $c_1$ with $q_{1n}^1$ as hyperedge feature; the $r_n^m x_n$ term in constraint $\delta_m$ is represented as a hyperedge covering vertices $v_1$, $v_2$ and $c_m$ with feature $r_n^m$.

\subsection{UniEGNN as an Interior-Point Method}
\label{sec:UniEGNN_as_Interior_Point_Method}

\begin{figure}[t]
    \centering
    \begin{minipage}{0.68\linewidth}
        \includegraphics[width=\linewidth]{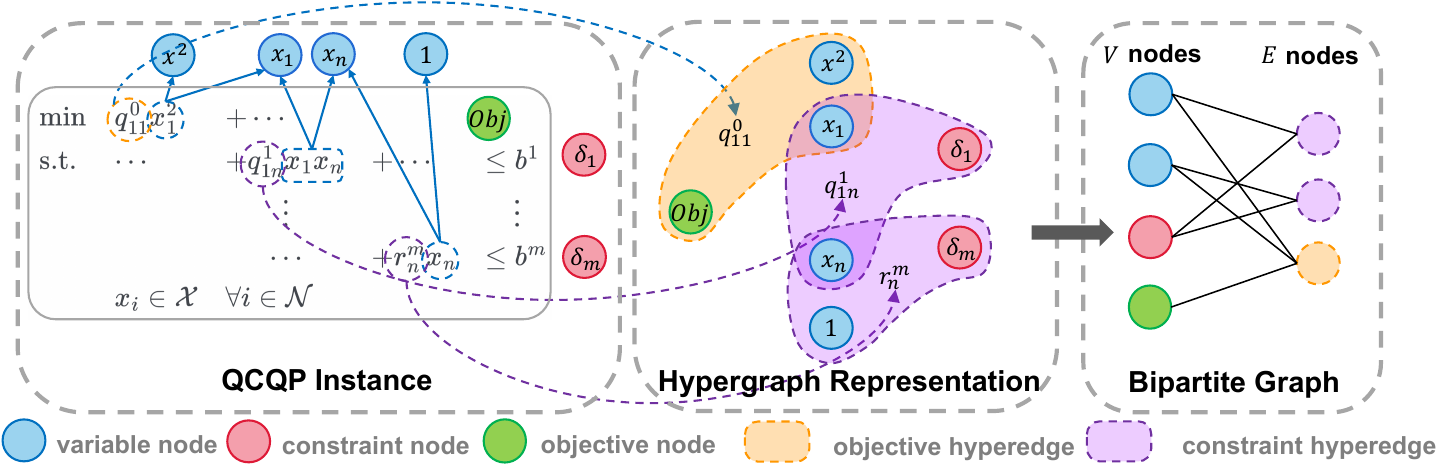}
        \caption{Variable relational hypergraph representation.}
        \label{fig:hypergraph_representation}
    \end{minipage}
    \hfill
    \begin{minipage}{0.29\linewidth}
        \centering
        \includegraphics[width=\linewidth]{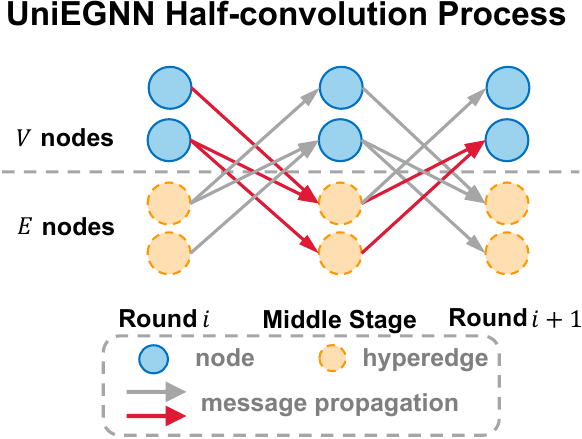}
        \caption{Half-convolution.}
        \label{fig:UniEGNN_convolution}
    \end{minipage}
\end{figure}

We then propose UniEGNN, which employs features of both vertices and hyperedges, converts the hypergraph into a bipartite graph, and updates vertex embeddings via spatial convolution (Section \ref{sec:hypergraph_basics}). Initially, different dimensional features of vertices and hyperedges are mapped into the same high-dimensional hidden space through a Multi-Layer Perceptron (MLP) layer. Then, vertices and hyperedges are reformulated as a bipartite graph via star expansion so that hyperedge features can be fully utilized in spatial convolution. The bipartite graph contains two sets of nodes, \( V \) nodes, and \( E \) nodes, which represent the original nodes and hyperedges, respectively. All edges connected with one \( E \) node share the same feature, as illustrated in Figure \ref{fig:hypergraph_representation}. Finally, given the bipartite graph reformulation, a half-convolution strategy is employed, as detailed in Equation \ref{eq:UniEGNN}:

\begin{equation}
    \begin{aligned}
        \text { (UniEGNN) }\left\{\begin{array}{l}
        h_e^{(t)}=\phi_e\left(h_e^{(t-1)},\psi_e\left(\{h_\alpha^{(t-1)}\}_{\alpha \in \mathcal{N}_v(e)}\right)\right)\\
        h_v^{(t)}=\phi_v\left(h_v^{(t-1)},\psi_v\left(\{h_\beta^{(t)}\}_{\beta \in \mathcal{N}_e(v)}\right)\right)
        \end{array}\right.,
        \label{eq:UniEGNN}
    \end{aligned}
\end{equation}

where \( h_e^{(t)} \) and \( h_v^{(t)} \) are hyperedge features and vertex features after the \( t \)-th convolution, respectively; \( \phi_e \) and \( \phi_v \) are permutation-invariant functions implemented by MLP; \( \phi_e \) and \( \psi_v \) are aggregation functions for hyperedges and vertices, which are \(\operatorname{SUM}\) and \(\operatorname{MEAN}\) in NeuralQP, respectively; \( \mathcal{N}_v(e) \) and \( \mathcal{N}_e(v) \) are the neighborhood of hyperedge \( e \) and vertex \( v \) defined in Section \ref{sec:hypergraph_basics}. Figure \ref{fig:UniEGNN_convolution} illustrates the message propagation in the half-convolution process, where features of both \( E \) nodes and \( V \) nodes are updated. The red arrow indicates the flow to the specific node.

\citet{qian2024exploring} proved the equivalence of MPNNs and IPMs for LPs by demonstrating that each step in the IPM can be replicated by MPNNs. We extend their proof to QPs in a similar manner. Firstly, IPM methods can efficiently solve QPs, and we provide a practical algorithm in Algorithm \ref{alg:IPM}. Secondly, to efficiently solve Equation \ref{eq:QP-reduced-KKT}in Algorithm \ref{alg:IPM}, we employ a conjugate gradient method (Algorithm \ref{alg:CG}). Then, to prove the equivalence of MPNNs to IPMs, we show that each step of the IPM can be replicated by MPNNs, following the approach of \citet{qian2024exploring}. Specifically, we first demonstrate that there exists an MPNN, $f_{\text{CG}}$, that replicates the conjugate gradient method (Lemma \ref{lemma:CG}). We then show that there exists an MPNN, $f_{\text{IPM}}$, that replicates the IPM (Theorem \ref{theorem:UniEGNN_IPM}).

\begin{lemma}
    There exists a MPNN \( f_{\text{CG}} \) composed of a finite number of message-passing steps that reproduces the conjugate gradient method described in Algorithm \ref{alg:CG}. Specifically, for any QP instance \( I = (\mQ, \vc, \mA, \vb) \) and any initial point \((\vx, \vy, \vz, \vw) > 0\), \( f_{\text{CG}} \) maps the graph \( \mathcal{H}(I) \) carrying the initial values on the variable and constraint nodes to the graph \( \mathcal{H}(I) \) carrying the output \(\Delta \vy\) of Algorithm \ref{alg:CG} on the variable nodes.
    \label{lemma:CG}
\end{lemma}

\begin{theorem}
    There exists an MPNN \( f_{\text{IPM}} \) composed of a finite number of message-passing steps that reproduces each iteration of the IPM described in Algorithm \ref{alg:IPM}. Specifically, for any QP instance \( I = (\mQ, \vc, \mA, \vb) \) and any iteration step \( t \geq 0 \), \( f_{\text{IPM}} \) maps the hypergraph \( \mathcal{H}(I) \) carrying the current iterate values on the variable and constraint nodes to the hypergraph \( \mathcal{H}(I) \) carrying the next iterate values \((\vx_{t+1}, \vz_{t+1}), (\vy_{t+1},\vw_{t+1})\) on the variable and constraint nodes.
    \label{theorem:UniEGNN_IPM}
\end{theorem}

Given that UniEGNN is a specific instance of MPNNs, the equivalence of MPNNs to IPMs directly implies the equivalence of UniEGNN to IPMs. Detailed proofs are provided in Appendix \ref{appendix:proof_equivalence_of_MPNNs_and_IPMs}.

\section{The General Hypergraph-based Optimization Framework}
\label{sec:general_hypergraph_based_optimization_framework}

\begin{figure}[t]
    \centering
    \includegraphics[width=1.0\linewidth]{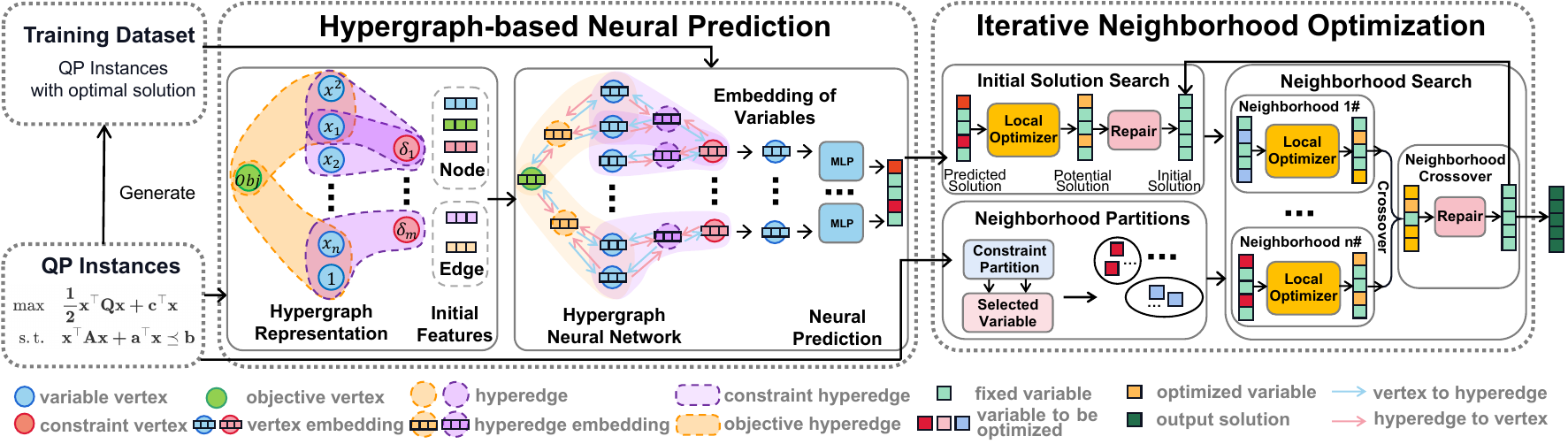}
    \caption{An overview of the NeuralQP framework. During hypergraph-based Neural Prediction, the QCQP is encoded into a variable relational hypergraph with initial vertex and hyperedge embeddings from the original problem. Then UniEGNN generates neural embeddings for each variable, utilizing both vertex and hyperedge features by converting the hypergraph into a bipartite graph. An MLP layer then predicts the optimal solution based on these embeddings. In Parallel Neighborhood Optimization, predicted solutions are first relaxed and repaired to obtain a feasible solution, followed by neighborhood partition, Parallel Neighborhood Optimization, and McCormick relaxation-based Q-Repair. The neighborhood search solution is again used as a feasible solution if the time limit is not reached; otherwise, the incumbent solution is output as the optimization result.}
    \label{fig:framework}
\end{figure}

This section presents our general hypergraph-based optimization framework for large-scale QCQPs, named NeuralQP. The framework of NeuralQP (illustrated in Figure \ref{fig:framework}) consists of \textit{Neural Prediction} (Section \ref{sec:neural_prediction}) which predicts a high-quality initial solution and \textit{Parallel Neighborhood Optimization} (Section \ref{sec:parallel_neighborhood_optimization}) which optimizes the incumbent solution by neighborhood search and crossover. 

\subsection{Neural Prediction}
\label{sec:neural_prediction}
During this stage, QCQPs are initially converted into variable relational hypergraphs (Section \ref{sec:hypergraph_based_representation_for_qcqp}). During the training process, \textit{UniEGNN} (Section \ref{sec:UniEGNN_as_Interior_Point_Method}) is applied on multiple hypergraphs to learn the typical structures of these problems. The outputs are interpreted as a probability of the optimal solution values similar to Neural Diving \citep{nairSolvingMixedInteger2021}, which are adopted as a heuristic for \textit{Parallel Neighborhood Optimization} in Section \ref{sec:parallel_neighborhood_optimization}. Details are presented in Appendix \ref{appendix:neural_prediction}.

\subsection{Parallel Neighborhood Optimization}
\label{sec:parallel_neighborhood_optimization}
In this stage, a variable proportion $\alpha_{ub} \in (0,1)$ is defined so that a large-scale QCQP with $n$ decision variables can be solved by a small-scale optimizer with $\alpha_{ub} n$ decision variables. Leveraging predicted values based on hypergraph neural networks, an initial feasible solution is searched by dynamically adjusting the radius using the Q-Repair algorithm (Section \ref{sec:initial_solution}). Then, a neighborhood search with a fixed radius is performed with an adaptive neighborhood partition strategy (Section \ref{sec:neighborhood-patition}). Finally, for multiple sets of neighborhood search solutions, neighborhood crossover and the Q-Repair algorithm are used to repair until the time limit or objective limit is reached. Finally, the current solution is output as the final optimization result.

\subsubsection{Q-Repair-Based Initial Feasible Solution}
\label{sec:initial_solution}
For a QCQP involving $n$ decision variables, NeuralQP arranges these variables in ascending order based on their predicted losses (Appendix \ref{appendix:model_formulation}). Let the initial search proportion be $\alpha \in (0,1)$. The first $(1-\alpha)n$ variables are fixed at their current values and the remaining variables are selected for optimization. The Q-Repair algorithm is then applied to identify constraints that are expected to be violated, and the neighborhood search radius is increased. The number of unfixed variables after the expansion is ensured to remain below or equal to $\alpha_{ub} n$. Therefore, a small-scale optimizer can be used to effectively find the initial feasible solution.

\begin{equation}
\label{eq:qrepair}
    \scalebox{0.85}{$
    \begin{aligned}
        \min\quad & \vx^{\mathrm{T}} \mQ^0 \vx+\left(\vr^0\right)^{\mathrm{T}} \vx \\
        \text {s.t.}\quad & \vx^{\mathrm{T}} \mQ^i \vx+\left(\vr^i\right)^{\mathrm{T}} \vx \le b_i\\
        & l_j \leq x_j \leq u_j \\
        & x_j = \hat{x}_j, \quad \forall x_j \in \mathcal{F}\\
        & \forall i \in \mathcal{M}, \quad \forall j \in \mathcal{N}
    \end{aligned}
    \begin{aligned}
        \xrightarrow[u_{jk},l_{jk} = \operatorname{McCor}(u_j,u_k,l_j,l_k)]{\bm{\phi}_{jk}:=x_jx_k,\forall j,k \in \mathcal{N}}
    \end{aligned}
    \begin{aligned}
        \min\quad & \vx^{\mathrm{T}} \mQ^0 \vx+\left(\vr^0\right)^{\mathrm{T}} \vx \\
        \text {s.t.}\quad & \sum q^i_{jk}\phi_{jk}+\left(\vr^i\right)^{\mathrm{T}} \vx \le b_i\\
        & l_j \leq x_j \leq u_j,l_{jk}\leq \phi_{jk}\leq u_{jk}\\
        & x_j = \hat{x}_j,\phi_{jk}=\hat{x}_jx_k,\forall x_j\in \mathcal{F}\\
        & \phi_{jk}=\hat{x}_j\hat{x}_k, \quad \forall x_j,x_k \in \mathcal{F}\\
        & \forall i \in \mathcal{M}, \quad \forall j,k \in \mathcal{N}.
    \end{aligned}
    $}
\end{equation}

The Q-Repair algorithm is based on the term-wise McCormick relaxation approach (Section \ref{sec:McCormick_relaxation}). 
First, the quadratic terms within the constraints are linearized through McCormick relaxation, which is shown in Equation \ref{eq:qrepair}, where $\hat{x}_j$ denotes the current solution of $x_i$ and $\mathcal{F}$ denotes the set of fixed variables. Then the linear repair algorithm \citep{ye2023gnn} is employed. For unfixed variables, their bounds are determined based on their original problem bounds. For fixed variables, their current solution values are utilized as the bounds for the respective terms.
% For unfixed variables, their bounds are used as the bounds for that term. As for fixed variables, their current solution is used as the bounds for that term. 
These bounds are all summed together to calculate the bounds for the left side of the constraint and be compared with the right-side coefficient. If the constraint is violated, the fixed variables within it are added to the neighborhood one by one until the solution becomes feasible. The related pseudocode is shown in Appendix \ref{appendix:Q-Repair}.

\subsubsection{Neighborhood Partition}
\label{sec:neighborhood-patition}
To reduce the number of violated constraints, neighborhood partitioning is carried out using the Adaptive Constraint Partition (ACP) framework \citep{ye2023adaptive}.
After randomly shuffling constraints, NeuralQP adds each variable of a constraint to a particular neighborhood one by one until the size of that neighborhood reaches its upper limit.
Due to the fixed radius of the neighborhoods, the number of neighborhoods is adaptive. More details and the pseudocode are shown in Appendix \ref{appendix:neighborhood_partition}.

The ACP framework increases the likelihood that variables within the same constraint are assigned to the same neighborhood, reducing the probability of constraints being violated. However, in QCQPs with a large number of variables within a single constraint, this approach may result in an excessive number of neighborhoods, which can reduce solving speed. Therefore, when the variable density (the average number of variables in one constraint) exceeds a certain threshold, a random neighborhood partitioning strategy is employed. All variables are randomly shuffled and neighborhood partitioning is performed sequentially so that each variable appears in only one neighborhood.

This adaptive neighborhood partitioning strategy, on one hand, retains the advantages of the ACP framework, ensuring the feasibility of subproblems after neighborhood partitioning. On the other hand, it addresses the issue of excessive neighborhood partitioning in dense problems within the ACP framework, thus accelerating the solving speed.

\subsubsection{Neighborhood Search and Crossover}
\label{sec:search-and-crossover}
After obtaining the initial feasible solution and neighborhood partitioning, NeuralQP uses small-scale solvers to perform fixed-radius neighborhood searches within multiple neighborhoods in parallel. Additionally, the neighborhood crossover strategy has shown promising results in linear problems \citep{ye2023gnn}. Therefore, neighborhood crossovers are also proposed for QCQPs to prevent getting trapped in local optima caused by limitations of the search radius.

To be specific, after the neighborhood partitioning, $l$ subproblems are solved in parallel using a small-scale solver. 
Then, neighborhood crossovers are performed between every two subproblems, resulting in $\floor{l/2}$ neighborhoods to be explored. 
To address the potential infeasible issues caused by neighborhood crossovers, the Q-Repair algorithm in Section \ref{sec:initial_solution} is applied again to repair the neighborhoods. 
Finally, the $\floor{l/2}$ neighborhoods are solved in parallel using a small-scale solver, and the best among them is selected as the result of this round of neighborhood optimization. If the time limit is not reached, the process proceeds to the next round; otherwise, the incumbent solution is output as the final result. Related pseudocodes are shown in Appendix \ref{appendix:neighborhood_search_and_crossover}.

\section{Experiments}
\label{sec:experiments}

Experiments were performed on two benchmark problems and real-world library QPLIB \citep{furiniQPLIBLibraryQuadratic2019}. The two benchmark problems are Quadratic Multiple Knapsack Problem (QMKP)(\citep{kellerer2004multidimensional}) in Appendix \ref{appendix:QMKP} and Random Quadratically Constrained Quadratic Program (RandQCP) in Appendix \ref{appendix:RandQCP}. For QMKP and RandQCP, we generated problems of Mini, 1000, 2000, 5000, and 10000 scales. Details about dataset settings are provided in Appendix \ref{appendix:experiments_details}. Results of experiments on QPLIB are presented in Appendix \ref{appendix:QPLIB} and ablation studies are presented in Appendix \ref{appendix:ablation_studies}.

In the testing stage, we used SCIP and Gurobi as the baselines. To verify the generalization of NeuralQP, we trained on three 1000er problem scales and then tested on the same or larger scales. Then our model was employed to obtain initial solutions, followed by iterative neighborhood search where the small-scale solvers were restricted to $30\%$ and $50\%$ of the number of variables. To ensure fair comparison, experiments were conducted in the following two aspects following \citet{ye2023gnn}. To study the effectiveness of NeuralQP, the objective value was compared with Gurobi and SCIP within the same wall-clock time (Section \ref{sec:comparison_objective_value}); to verify the efficiency of NeuralQP, the running time till reaching the same objective value was compared with Gurobi and SCIP (Section \ref{sec:comparison_runtime}). To avoid the randomness of the experiment, we tested ten instances of each benchmark problem at every scale. Each instance was run five times, and the mean and standard deviation were recorded. More details of experiment settings are listed in the Appendix \ref{appendix:experiments_details}. The code and data for reproducing the experiment results have been uploaded to \href{https://anonymous.4open.science/r/NeuralQP-Anonymous-7243}{this repository}.

\subsection{Comparison of Objective Value}
\label{sec:comparison_objective_value}

This section compares our framework with large-scale solvers Gurobi and SCIP, evaluating their performance within the same runtime. As shown in Table \ref{table:1}, using a solver with a scale of only $30\% $of the number of variables, out method can obtain better objective value compared to baselines on all benchmark problems except QMKP of 1000 scale. When the restriction increases to $50\%$, our method typically performs better. For small-scale problems, the advantages of our framework compared to large-scale solvers are not pronounced; however, when it comes to solving large-scale problems, our framework significantly outperforms both large-scale solvers within the same time limit.

Furthermore, the problem scale in the training set does not have a significant impact, demonstrating the generalization ability of our method. Even if only the solutions to small-scale problems are known, our NeuralQP can still be used to solve large-scale problems. All the above analyses demonstrate that our method can use small-scale solvers to effectively address large-scale QCQPs.

\begin{table}[h]
\centering
\scriptsize
\caption{Comparison of objective values with SCIP and Gurobi within the same runtime. Ours-30\%S and Ours-30\%G mean the scale-limited versions of SCIP and Gurobi respectively, with the variable proportion $\alpha$ limited to 30\%. “↑” means the result is better than or equal to the baseline. Each value is averaged among 10 instances.}
\label{table:1}
\resizebox{\textwidth}{!}{%
\begin{tabular}{cccccccccc}
\hline
                            \multirow{2}{*}{Method}     & \multirow{2}{*}{Train set}& \multicolumn{4}{c}{QMKP}      & \multicolumn{4}{c}{RandQCP}   \\ \cline{3-10} 
                            &        & 1000 & 2000 & 5000 & 10000 & 1000 & 2000 & 5000 & 10000 \\ \hline
SCIP                        &        & \textbf{1854.65} & 2699.15 & 8880.72 & 11336.42 & 186.27 & 343.41 & 806.12 & 1408.09 \\
\multirow{3}{*}{Ours-30\%S} & Mini   & 1647.08 & 2816.87↑ & 13781.18↑ & 14051.19↑ & 295.81↑ & 579.26↑ & 1494.55↑ & 2984.41↑ \\
                            & 1000  & 1647.61 & 2766.66↑ & 13771.44↑ & 14050.79↑ & 295.48↑ & 578.35↑ & 1495.09↑ & 2984.73↑ \\
                            & 2000 & - & 2815.14↑ & 13782.54↑ & 14051.37↑ & - & 580.24↑ & 1494.35↑ & 2982.67↑ \\
\multirow{3}{*}{Ours-50\%S} & Mini   & 1711.71 & \textbf{2843.69↑} & 13795.12↑ & 14066.36↑ & \textbf{297.60↑} & 581.60↑ & \textbf{1499.58↑} & \textbf{2993.30↑} \\
                            & 1000  & 1712.77 & 2818.41↑ & 13793.15↑ & 14067.14↑ & 297.21↑ & 581.30↑ & 1499.56↑ & 2992.54↑ \\
                            & 2000 & - & 2843.47↑ & \textbf{13796.34↑} & \textbf{14066.81↑} & - & \textbf{582.50↑} & 1498.82↑ & 2992.68↑ \\ \hline
Gurobi                      &      & 1499.99 & 2634.47 & 9706.71 & 11199.70 & 285.65 & 554.21 & 1435.00 & 2807.03 \\
\multirow{3}{*}{Ours-30\%G} & Mini   & \textbf{1838.82↑} & 2816.60↑ & 13783.47↑ & 14045.67↑ & 296.55↑ & 580.86↑ & 1495.28↑ & 2985.00↑ \\
                            & 1000  & 1802.04↑ & 2768.76↑ & 13775.64↑ & 14052.28↑ & 295.27↑ & 581.07↑ & 1495.11↑ & 2983.48↑ \\
                            & 2000 & - & 2817.28↑ & 13782.34↑ & 14052.23↑ & - & 581.12↑ & 1496.26↑ & 2983.09↑ \\
\multirow{3}{*}{Ours-50\%G} & Mini   & 1778.29↑ & 2844.58↑ & \textbf{13797.71↑} & 14061.14↑ & \textbf{298.35↑} & \textbf{583.17↑} & 1499.81↑ & \textbf{2992.92↑} \\
                            & 1000  & 1751.77↑ & 2825.65↑ & 13793.91↑ & \textbf{14068.15↑} & 297.15↑ & 582.47↑ & \textbf{1500.12↑} & 2992.66↑ \\
                            & 2000 & - & \textbf{2844.46↑} & 13797.12↑ & 14067.17↑ & - & 582.73↑ & 1500.06↑ & 2992.49↑ \\ \hline
Time                        &      & 100s & 600s & 1800s & 3600s & 100s & 600s & 1800s & 3600s \\ \hline
\end{tabular}
}
\end{table}

\subsection{Comparison of Runtime}
\label{sec:comparison_runtime}

In addition to solving performance, we also evaluated the time efficiency of NeuralQP by comparing the duration required to reach equivalent optimization results. According to the results in Table \ref{table:2}, NeuralQP consistently requires less time to achieve the same objective values compared to Gurobi and SCIP across all models and all benchmark problems, except for QMKP of 1000-scale. Even though the objective values achieved within a similar time limit may appear close, NeuralQP typically achieves the same performance as the baselines with only 10\% or even 1\% the time required.

\begin{table}[h]
\centering
\scriptsize
\caption{Comparsion of running time with SCIP and Gurobi till the same objective value. Notations are similar to Table \ref{table:1}. $>$ indicates the inability to achieve the target objective function within the given time limit for some instances (Appendix \ref{appendix:experiment_settings}). Each value is averaged among 10 instances.}
\resizebox{\textwidth}{!}{
\label{table:2}
\begin{tabular}{cccccccccc}
\hline
                            \multirow{2}{*}{Method}     & \multirow{2}{*}{Train set}& \multicolumn{4}{c}{QMKP}      & \multicolumn{4}{c}{RandQCP}   \\ \cline{3-10} 
                            &        & 1000 & 2000 & 5000 & 10000 & 1000 & 2000 & 5000 & 10000 \\ \hline
SCIP                        &       & \textbf{27.44s} & $>$2920s & $>$10804s & $>$21600s & $>$965s & $>$3600s & $>$10800s & $>$21600s \\
\multirow{3}{*}{Ours-30\%S} & Mini   & 122.03s & 397.85s & 891.56s & 1790.94s & 36.18s & 119.91s & 446.19s & 847.82s\\
                            & 1000  &  121.70s & 574.73s & 1147.64s & 1821.29s & 38.45s & 120.87s & 489.88s & 846.42s\\
                            & 2000 & - & 393.61s & 871.83 & 1809.83s & - & 106.21s & 393.76s & 960.09s \\
\multirow{3}{*}{Ours-50\%S} & Mini   &  97.54s & \textbf{278.14s} & 637.21s & 1170.38s & 32.51s & 57.40s & \textbf{157.67s} & \textbf{324.01s} \\
                            & 1000  & 91.24s & 434.20s & 801.97s & 1160.33s & \textbf{24.74s} & \textbf{51.90s} & 227.88s & 377.23s \\
                            & 2000 & - & 290.88s & \textbf{608.08s} & \textbf{1160.13s} & - & 53.87s & 188.60s & 340.88s \\ \hline
Objective                   &      & 1604.78 & 2740.92 & 13677.24 & 13986.67 & 289.67 & 565.40 & 1479.21 & 2947.39 \\ \hline
Gurobi                      &        & 583.30s & 1687.50s & $>$10786s & $>$17206s & 521.60s & 2423.70s & 9588.00s & 19813.80s  \\
\multirow{3}{*}{Ours-30\%G} & Mini   & 48.10s & 387.32s & 885.34s & 1602.41s & 61.29s & 154.84s & 259.04s & 447.55s \\
                            & 1000  & \textbf{45.36s} & 560.12s & 1093.89s & 1604.20s & 54.95s & 434.94s & 203.81s & 454.21s \\
                            & 2000 & - & 394.61s & 858.90s & 1599.50s & - & 316.79s & 630.31s & 420.41s \\
\multirow{3}{*}{Ours-50\%G} & Mini   & 79.48s & \textbf{276.43s} & 628.16s & 1084.86s & \textbf{27.07s} & 100.36s & \textbf{203.03s} & 495.87s \\
                            & 1000  & 74.15s & 390.11s & 794.12s & 1074.31s & 35.98s & 127.03s & 252.42s & \textbf{158.53s} \\
                            & 2000 & - & 283.94s & \textbf{599.56s} & \textbf{1067.82s} & - & \textbf{88.30s} & 304.49s & 174.84s \\ \hline
Objective                   &      & 1689.39 & 2738.92 & 13676.15 & 13974.16 & 290.72 & 567.74 & 1477.62 & 2944.70 \\ \hline
\end{tabular}
}
\end{table}

\section{Conclusion}
This paper presents NeuralQP, a pioneering hypergraph-based optimization framework for large-scale QCQPs. Key features of NeuralQP include: 1) a variable relational hypergraph as a complete representation of QCQPs without any assumption and 2) a McCormick relaxation-based repair algorithm that can identify illegal constraints. Experimental results demonstrate that NeuralQP achieves equivalent quality solutions in less than 10\% of the time compared to leading solvers on large-scale QCQPs, setting the groundwork for solving nonlinear problems via machine learning. Despite these advances, the framework is limited by the lack of publicly available datasets. Future work should focus on expanding dataset availability, extending the proof to nonlinear optimization, and generalizing to problems with integer variables. Please refer to Appendix \ref{appendix:limitation_impacts_future_directions} for detailed analysis.

\bibliography{reference}

\begin{thebibliography}{37}
\providecommand{\natexlab}[1]{#1}
\providecommand{\url}[1]{\texttt{#1}}
\expandafter\ifx\csname urlstyle\endcsname\relax
  \providecommand{\doi}[1]{doi: #1}\else
  \providecommand{\doi}{doi: \begingroup \urlstyle{rm}\Url}\fi

\bibitem[Annamalai(2016)]{annamalai2016finding}
C.~Annamalai.
\newblock Finding perfect matchings in bipartite hypergraphs.
\newblock In \emph{Proceedings of the twenty-seventh annual ACM-SIAM symposium on Discrete algorithms}, pages 1814--1823. SIAM, 2016.

\bibitem[Balas(1969)]{balas1969duality}
E.~Balas.
\newblock Duality in discrete programming: Ii. the quadratic case.
\newblock \emph{Management Science}, 16\penalty0 (1):\penalty0 14--32, 1969.

\bibitem[Bertsimas and Stellato(2020)]{bertsimasVoiceOptimization2020}
D.~Bertsimas and B.~Stellato.
\newblock The {{Voice}} of {{Optimization}}, June 2020.

\bibitem[Bertsimas and Stellato(2021)]{bertsimasOnlineMixedIntegerOptimization2021}
D.~Bertsimas and B.~Stellato.
\newblock Online {{Mixed-Integer Optimization}} in {{Milliseconds}}, Mar. 2021.

\bibitem[Bonami et~al.(2018)Bonami, Lodi, and Zarpellon]{bonamiLearningClassificationMixedInteger2018}
P.~Bonami, A.~Lodi, and G.~Zarpellon.
\newblock Learning a {{Classification}} of {{Mixed-Integer Quadratic Programming Problems}}.
\newblock In W.-J. {van Hoeve}, editor, \emph{Integration of {{Constraint Programming}}, {{Artificial Intelligence}}, and {{Operations Research}}}, Lecture {{Notes}} in {{Computer Science}}, pages 595--604, {Cham}, 2018. {Springer International Publishing}.

\bibitem[Bonami et~al.(2022)Bonami, Lodi, and Zarpellon]{bonamiClassifierDecideLinearization2022}
P.~Bonami, A.~Lodi, and G.~Zarpellon.
\newblock A {{Classifier}} to {{Decide}} on the {{Linearization}} of {{Mixed-Integer Quadratic Problems}} in {{CPLEX}}.
\newblock \emph{Operations Research}, 70\penalty0 (6):\penalty0 3303--3320, Nov. 2022.

\bibitem[Bretto(2013)]{bretto2013hypergraph}
A.~Bretto.
\newblock Hypergraph theory.
\newblock \emph{An introduction. Mathematical Engineering. Cham: Springer}, 1, 2013.

\bibitem[Burer and Letchford(2012)]{burer2012non}
S.~Burer and A.~N. Letchford.
\newblock Non-convex mixed-integer nonlinear programming: A survey.
\newblock \emph{Surveys in Operations Research and Management Science}, 17\penalty0 (2):\penalty0 97--106, 2012.

\bibitem[Burer and Saxena(2012)]{burerMILPRoadMIQCP2012}
S.~Burer and A.~Saxena.
\newblock The {{MILP Road}} to {{MIQCP}}.
\newblock In J.~Lee and S.~Leyffer, editors, \emph{Mixed {{Integer Nonlinear Programming}}}, The {{IMA Volumes}} in {{Mathematics}} and Its {{Applications}}, pages 373--405, {New York, NY}, 2012. {Springer}.
\newblock ISBN 978-1-4614-1927-3.
\newblock \doi{10.1007/978-1-4614-1927-3_13}.

\bibitem[Burkard et~al.(1997)Burkard, Karisch, and Rendl]{burkard1997qaplib}
R.~E. Burkard, S.~E. Karisch, and F.~Rendl.
\newblock Qaplib--a quadratic assignment problem library.
\newblock \emph{Journal of Global optimization}, 10:\penalty0 391--403, 1997.

\bibitem[Chen et~al.(2023)Chen, Liu, Wang, Lu, and Yin]{chenRepresentingMixedIntegerLinear2023}
Z.~Chen, J.~Liu, X.~Wang, J.~Lu, and W.~Yin.
\newblock On {{Representing Mixed-Integer Linear Programs}} by {{Graph Neural Networks}}, May 2023.

\bibitem[Coxeter(1950)]{coxeter1950self}
H.~S. Coxeter.
\newblock Self-dual configurations and regular graphs.
\newblock \emph{Bulletin of the American Mathematical Society}, 56\penalty0 (5):\penalty0 413--455, 1950.

\bibitem[Dai and Gao(2023)]{daiHypergraphComputation2023}
Q.~Dai and Y.~Gao.
\newblock \emph{Hypergraph {{Computation}}}.
\newblock Artificial {{Intelligence}}: {{Foundations}}, {{Theory}}, and {{Algorithms}}. {Springer Nature Singapore}, {Singapore}, 2023.

\bibitem[Ding et~al.(2019)Ding, Zhang, Shen, Li, Wang, Xu, and Song]{dingAcceleratingPrimalSolution2019}
J.-Y. Ding, C.~Zhang, L.~Shen, S.~Li, B.~Wang, Y.~Xu, and L.~Song.
\newblock Accelerating {{Primal Solution Findings}} for {{Mixed Integer Programs Based}} on {{Solution Prediction}}, Sept. 2019.

\bibitem[Elloumi and Lambert(2019)]{elloumi2019global}
S.~Elloumi and A.~Lambert.
\newblock Global solution of non-convex quadratically constrained quadratic programs.
\newblock \emph{Optimization methods and software}, 34\penalty0 (1):\penalty0 98--114, 2019.

\bibitem[Furini et~al.(2019)Furini, Traversi, Belotti, Frangioni, Gleixner, Gould, Liberti, Lodi, Misener, Mittelmann, Sahinidis, Vigerske, and Wiegele]{furiniQPLIBLibraryQuadratic2019}
F.~Furini, E.~Traversi, P.~Belotti, A.~Frangioni, A.~Gleixner, N.~Gould, L.~Liberti, A.~Lodi, R.~Misener, H.~Mittelmann, N.~V. Sahinidis, S.~Vigerske, and A.~Wiegele.
\newblock {{QPLIB}}: A library of quadratic programming instances.
\newblock \emph{Mathematical Programming Computation}, 11\penalty0 (2):\penalty0 237--265, June 2019.

\bibitem[Gallo et~al.(1980)Gallo, Hammer, and Simeone]{gallo1980quadratic}
G.~Gallo, P.~L. Hammer, and B.~Simeone.
\newblock Quadratic knapsack problems.
\newblock \emph{Combinatorial optimization}, pages 132--149, 1980.

\bibitem[Galloway et~al.(2015)Galloway, Sreenath, Ames, and Grizzle]{galloway2015torque}
K.~Galloway, K.~Sreenath, A.~D. Ames, and J.~W. Grizzle.
\newblock Torque saturation in bipedal robotic walking through control lyapunov function-based quadratic programs.
\newblock \emph{IEEE Access}, 3:\penalty0 323--332, 2015.

\bibitem[Gao et~al.(2023)Gao, Feng, Ji, and Ji]{gaoHGNNGeneralHypergraph2023}
Y.~Gao, Y.~Feng, S.~Ji, and R.~Ji.
\newblock {{HGNN}}+: {{General Hypergraph Neural Networks}}.
\newblock \emph{IEEE Transactions on Pattern Analysis and Machine Intelligence}, 45\penalty0 (3):\penalty0 3181--3199, Mar. 2023.

\bibitem[Gasse et~al.(2019)Gasse, Ch{\'e}telat, Ferroni, Charlin, and Lodi]{gasseExactCombinatorialOptimization2019}
M.~Gasse, D.~Ch{\'e}telat, N.~Ferroni, L.~Charlin, and A.~Lodi.
\newblock Exact {{Combinatorial Optimization}} with {{Graph Convolutional Neural Networks}}, Oct. 2019.

\bibitem[Ghaddar et~al.(2022)Ghaddar, {G{\'o}mez-Casares}, {Gonz{\'a}lez-D{\'\i}az}, {Gonz{\'a}lez-Rodr{\'\i}guez}, {Pateiro-L{\'o}pez}, and {Rodr{\'\i}guez-Ballesteros}]{ghaddarLearningSpatialBranching2022}
B.~Ghaddar, I.~{G{\'o}mez-Casares}, J.~{Gonz{\'a}lez-D{\'\i}az}, B.~{Gonz{\'a}lez-Rodr{\'\i}guez}, B.~{Pateiro-L{\'o}pez}, and S.~{Rodr{\'\i}guez-Ballesteros}.
\newblock Learning for {{Spatial Branching}}: {{An Algorithm Selection Approach}}, Apr. 2022.

\bibitem[Gondzio and Grothey(2007)]{gondzio2007parallel}
J.~Gondzio and A.~Grothey.
\newblock Parallel interior-point solver for structured quadratic programs: Application to financial planning problems.
\newblock \emph{Annals of Operations Research}, 152:\penalty0 319--339, 2007.

\bibitem[Hiley and Julstrom(2006)]{hiley2006quadratic}
A.~Hiley and B.~A. Julstrom.
\newblock The quadratic multiple knapsack problem and three heuristic approaches to it.
\newblock In \emph{Proceedings of the 8th annual conference on Genetic and evolutionary computation}, pages 547--552, 2006.

\bibitem[Huang and Yang(2021)]{huangUniGNNUnifiedFramework2021}
J.~Huang and J.~Yang.
\newblock {{UniGNN}}: A {{Unified Framework}} for {{Graph}} and {{Hypergraph Neural Networks}}, May 2021.

\bibitem[Ichnowski et~al.(2021)Ichnowski, Jain, Stellato, Banjac, Luo, Borrelli, Gonzalez, Stoica, and Goldberg]{ichnowski2021accelerating}
J.~Ichnowski, P.~Jain, B.~Stellato, G.~Banjac, M.~Luo, F.~Borrelli, J.~E. Gonzalez, I.~Stoica, and K.~Goldberg.
\newblock Accelerating quadratic optimization with reinforcement learning.
\newblock \emph{Advances in Neural Information Processing Systems}, 34:\penalty0 21043--21055, 2021.

\bibitem[Kannan et~al.(2023)Kannan, Nagarajan, and Deka]{kannanLearningAccelerateGlobal2023}
R.~Kannan, H.~Nagarajan, and D.~Deka.
\newblock Learning to {{Accelerate}} the {{Global Optimization}} of {{Quadratically-Constrained Quadratic Programs}}, Feb. 2023.

\bibitem[Kellerer et~al.(2004)Kellerer, Pferschy, Pisinger, Kellerer, Pferschy, and Pisinger]{kellerer2004multidimensional}
H.~Kellerer, U.~Pferschy, D.~Pisinger, H.~Kellerer, U.~Pferschy, and D.~Pisinger.
\newblock \emph{Multidimensional knapsack problems}.
\newblock Springer, 2004.

\bibitem[McCormick(1976)]{mccormickComputabilityGlobalSolutions1976}
G.~P. McCormick.
\newblock Computability of global solutions to factorable nonconvex programs: {{Part I}} \textemdash{} {{Convex}} underestimating problems.
\newblock \emph{Mathematical Programming}, 10\penalty0 (1):\penalty0 147--175, Dec. 1976.
\newblock ISSN 1436-4646.
\newblock \doi{10.1007/BF01580665}.

\bibitem[Nair et~al.(2021)Nair, Bartunov, Gimeno, {von Glehn}, Lichocki, Lobov, O'Donoghue, Sonnerat, Tjandraatmadja, Wang, Addanki, Hapuarachchi, Keck, Keeling, Kohli, Ktena, Li, Vinyals, and Zwols]{nairSolvingMixedInteger2021}
V.~Nair, S.~Bartunov, F.~Gimeno, I.~{von Glehn}, P.~Lichocki, I.~Lobov, B.~O'Donoghue, N.~Sonnerat, C.~Tjandraatmadja, P.~Wang, R.~Addanki, T.~Hapuarachchi, T.~Keck, J.~Keeling, P.~Kohli, I.~Ktena, Y.~Li, O.~Vinyals, and Y.~Zwols.
\newblock Solving {{Mixed Integer Programs Using Neural Networks}}, July 2021.

\bibitem[Qian et~al.(2024)Qian, Ch{\'e}telat, and Morris]{qian2024exploring}
C.~Qian, D.~Ch{\'e}telat, and C.~Morris.
\newblock Exploring the power of graph neural networks in solving linear optimization problems.
\newblock In \emph{International Conference on Artificial Intelligence and Statistics}, pages 1432--1440. PMLR, 2024.

\bibitem[R{\"o}dl and Skokan(2004)]{rodl2004regularity}
V.~R{\"o}dl and J.~Skokan.
\newblock Regularity lemma for k-uniform hypergraphs.
\newblock \emph{Random Structures \& Algorithms}, 25\penalty0 (1):\penalty0 1--42, 2004.

\bibitem[Stellato et~al.(2017)Stellato, Banjac, Goulart, Bemporad, and Boyd]{stellatoOSQPOperatorSplitting2017}
B.~Stellato, G.~Banjac, P.~Goulart, A.~Bemporad, and S.~Boyd.
\newblock {{OSQP}}: {{An Operator Splitting Solver}} for {{Quadratic Programs}}, Nov. 2017.

\bibitem[Vanderbei(1999)]{vanderbeiLOQO1999}
R.~J. Vanderbei.
\newblock Loqo:an interior point code for quadratic programming.
\newblock \emph{Optimization Methods and Software}, 11\penalty0 (1-4):\penalty0 451--484, 1999.
\newblock \doi{10.1080/10556789908805759}.

\bibitem[Weisfeiler and Leman(1968)]{weisfeiler1968reduction}
B.~Weisfeiler and A.~Leman.
\newblock The reduction of a graph to canonical form and the algebra which appears therein.
\newblock \emph{nti, Series}, 2\penalty0 (9):\penalty0 12--16, 1968.

\bibitem[Ye et~al.(2023{\natexlab{a}})Ye, Wang, Xu, Wang, and Jiang]{ye2023adaptive}
H.~Ye, H.~Wang, H.~Xu, C.~Wang, and Y.~Jiang.
\newblock Adaptive constraint partition based optimization framework for large-scale integer linear programming (student abstract).
\newblock In \emph{Proceedings of the AAAI Conference on Artificial Intelligence}, volume~37, pages 16376--16377, 2023{\natexlab{a}}.

\bibitem[Ye et~al.(2023{\natexlab{b}})Ye, Xu, Wang, Wang, and Jiang]{ye2023gnn}
H.~Ye, H.~Xu, H.~Wang, C.~Wang, and Y.~Jiang.
\newblock Gnn\&gbdt-guided fast optimizing framework for large-scale integer programming.
\newblock 2023{\natexlab{b}}.

\bibitem[Zhang et~al.(2013)Zhang, Yao, Iu, Fernando, and Wong]{zhang2013sequential}
Y.~Zhang, F.~Yao, H.~H.-C. Iu, T.~Fernando, and K.~P. Wong.
\newblock Sequential quadratic programming particle swarm optimization for wind power system operations considering emissions.
\newblock \emph{Journal of Modern Power Systems and Clean Energy}, 1\penalty0 (3):\penalty0 227--236, 2013.

\end{thebibliography}

%%%%%%%%%%%%%%%%%%%%%%%%%%%%%%%%%%%%%%%%%%%%%%%%%%%%%%%%%%%%
% Appendix
%%%%%%%%%%%%%%%%%%%%%%%%%%%%%%%%%%%%%%%%%%%%%%%%%%%%%%%%%%%%
\newpage
\section*{Appendix}
\appendix

\section{Characteristics of QCQPs}
\label{appendix:quadratic_programming}
QCQPs encompass a broad class of optimization problems, where both the objective function and the constraints are quadratic. Equation \ref{eq:qp} lays down the foundational structure of QCQPs. In this section, we briefly introduce the characteristics of QCQPs, emphasizing both mathematical and machine learning-based approaches.

\subsection{Convexity in QCQPs}
\label{appendix:convexity}
Convexity plays a pivotal role in the realm of optimization, significantly impacting solvability and computational tractability. A QCQP is termed convex only if both the objective function and the feasible region are convex. In such scenarios, the global optimum can be efficiently found using polynomial-time algorithms. However, general QCQPs do not have the convexity property. QCQPs with integer variables are even harder, posing significant challenges in finding the optimal solution, which necessitates the use of heuristic methods, relaxation techniques, or local search strategies.

\subsection{Mixed-Integer and Binary Quadratic Programs}
\label{appendix:miqcp_bqcp}
When QCQPs incorporate integer constraints on certain variables, they become Mixed-Integer Quadratically Constrained Programs (MIQCPs). MIQCPs inherit the complexities of both integer programming and quadratic programming, which combine to make them particularly challenging to solve. Binary Quadratically Constrained Programs (BQCPs), a special case of MIQCPs, impose strict restrictions on variables by exclusively allowing binary values.

\subsection{Solution Methods}
\label{appendix:solution_strategy}
Solving QCQPs involves a variety of methods. We divide these methods into two main categories: \textit{mathematical approaches} and \textit{machine learning approaches}.

\textit{Mathematical approaches} to QCQPs typically involve convex optimization algorithms for convex instances and approximation methods, branch-and-bound techniques, or metaheuristic approaches for nonconvex instances. The choice of method depends heavily on the problem's structure, size, and the desired balance between solution quality and computational resources. For further reference, readers are encouraged to refer to key works in the field, such as \cite{burerMILPRoadMIQCP2012} and \citep{burer2012non}.

The use of \textit{machine learning methods} in solving QCQPs is notably scarce, with only a few studies, as is introduced in Section \ref{sec:introduction}. The problems in \cite{kannanLearningAccelerateGlobal2023} are relatively small compared to our experiments.  Other studies focus on Quadratic Programs (QPs) (i.e., problems with quadratic objective functions and linear constraints) \citep{bonamiLearningClassificationMixedInteger2018,bonamiClassifierDecideLinearization2022,ichnowski2021accelerating,bertsimasOnlineMixedIntegerOptimization2021} which are not directly comparable to our approach. Meanwhile, classic datasets like QPLIB \citep{furiniQPLIBLibraryQuadratic2019} exist but are limited in problem count (453 in total with 133 binary ones) and too diverse, making them unsuitable for training machine learning models due to the varied distribution of problems. Specialized datasets like QAPLIB \citep{burkard1997qaplib} also exist but focus on specific problem types.

In our method, to enable the use of small-scale solvers for general large-scale problems, we utilized two benchmark problems and adopted a random generation strategy. Such strategy ensures that 
\begin{enumerate}
    \item Our approach is tested and validated across different problem types;
    \item The problems generated are sufficiently large and difficult to solve.
\end{enumerate}

\section{Proof of the Equivalence of UniEGNN and IPMs}
\label{appendix:proof}

In this section, we demonstrate that MPNNs are equivalent to the general framework of the IPM for solving quadratic programming problems within our hypergraph representation. We first establish that there exists an MPNN \( f_{\text{CG}} \) that reproduces the conjugate gradient method outlined in Algorithm \ref{alg:CG} (Lemma \ref{lemma:CG}). Subsequently, we show that there exists an MPNN \( f_{\text{IPM}} \) that is equivalent to the IPM described in Algorithm \ref{alg:IPM} (Theorem \ref{theorem:UniEGNN_IPM}). Given that UniEGNN belongs to the class of MPNNs, we thus establish the equivalence between UniEGNN and IPM, validating its theoretical foundation.

\subsection{The General IPM for Nonlinear Optimization}

Consider a nonlinear optimization problem of the form:

\begin{equation}
    \begin{aligned}
        \min &\quad f(\vx) \\
        \text{subject to} &\quad h_i(\vx) \geq 0, \quad i = 1, \ldots, m,
    \end{aligned}
    \label{eq:def}
\end{equation}

where $\vx \in \mathbb{R}^n$, $f : \mathbb{R}^n \rightarrow \mathbb{R}$ and $h_i : \mathbb{R}^n \rightarrow \mathbb{R}$ are assumed to be twice continuously differentiable. Introduce slack variables $w_i$ to each of the constraints in Equation \ref{eq:def} and reformulate the problem as:

\begin{equation}
    \begin{aligned}
        \min &\quad f(\vx) - \mu \sum_{i=1}^{m} \log(w_i) \\
        \text{subject to} & \quad \vh(\vx) - \vw = \mathbf{0}, \quad \vw \geq \mathbf{0}
    \end{aligned}
\end{equation}

where $\vw$ is a vector of slack variables $w_i$. The Lagrangian for the problem is

\begin{equation}
    L(\vx, \vw, \vy, \mu) = f(\vx) - \mu \sum_{i=1}^{m} \log(w_i) - \vy^\top(\vh(\vx) - \vw),
\end{equation}

and the \textit{first-order KKT conditions} for a minimum are

\begin{equation}
    \begin{aligned}
        \nabla_x L &= \nabla f(\vx) - \nabla \vh(\vx)^\top \vy = \mathbf{0}, \\
        \nabla_w L &= -\mu \mW^{-1}\ve + \vy = \mathbf{0}, \\
        \nabla_y L &= \vh(\vx) - \vw = \mathbf{0},
        \label{eq:kkt}
    \end{aligned}
\end{equation}

where $\mW$ is the diagonal matrix with $w_i$ being the $i$-th diagonal element (i.e., $\mW = \operatorname{diag}(w_1, \dots, w_m)$) and $\ve$ is a vector of ones. Multiply the second equation of Equation \ref{eq:kkt} by $\mY$ and we get the \textit{primal-dual system}:

\begin{equation}
    \begin{aligned}
        \nabla f(\vx) - \nabla \vh(\vx)^\top \vy &= 0, \\
        -\mu \ve + \mW\mY \ve &= \mathbf{0}, \\
        \vh(\vx) - \vw &= \mathbf{0},
        \label{eq:primal-dual}
    \end{aligned}
\end{equation}

where again $\mY$ is constructed similarly to $\mW$ (we then do not mention it in the following). The basic numerical algorithm for finding a solution to the primal-dual system Equation \ref{eq:primal-dual} is to use \textit{Newton's method}. To simplify notations the following definitions are introduced:

\begin{equation*}
    \begin{aligned}
        \mH(\vx, \vy) &:= \nabla^2 f(\vx) - \sum_{i=1}^{m} y_i \nabla^2 h_i(\vx), \\
        \mA(\vx) &:= \nabla \vh(\vx).
    \end{aligned}
\end{equation*}

The Newton system for Equation (6) is

\begin{equation}
    \begin{bmatrix}
        \mH(\vx, \vy) & \mathbf{0} & -\mA(\vx)^\top \\
        \mathbf{0} & \mY & \mW \\
        \mA(\vx) & -\mI & \mathbf{0}
    \end{bmatrix}
    \begin{bmatrix}
        \Delta \vx \\
        \Delta \vw \\
        \Delta \vy
    \end{bmatrix}
    =
    \begin{bmatrix}
        -\nabla f(\vx) + \mA(\vx)^\top \vy \\
        \mu \ve - \mW\mY \ve \\
        -\vh(\vx) + \vw
    \end{bmatrix}.
    \label{eq:Newton-system}
\end{equation}

The system Equation \ref{eq:Newton-system} is not symmetric but is easily symmetrized by multiplying the second equation by $\mW^{-1}$ yielding
\begin{equation}
    \begin{bmatrix}
        -\mH(\vx, \vy) & \mathbf{0} & \mA(\vx)^\top \\
        \mathbf{0} & -\mW^{-1}\mY & -\mI \\
        \mA(\vx) & -\mI & \mathbf{0}
    \end{bmatrix}
    \begin{bmatrix}
        \Delta \vx \\
        \Delta \vw \\
        \Delta \vy
    \end{bmatrix}
    =
    \begin{bmatrix}
        \sigma \\
        -\gamma \\
        \rho
    \end{bmatrix},\text{ where } \begin{cases}
        \sigma &= \nabla f(\vx) - \mA(\vx)^\top \vy, \\
        \gamma &= \mu \mW^{-1}\ve - \vy, \\
        \rho &=  \vw - \vh(\vx).
    \end{cases}
\end{equation}

Here $\rho$ measures \textit{primal infeasibility} and $\sigma$ measures \textit{dual infeasibility} by analogy with linear programming. Since

\begin{equation*}
    \Delta \vw = \mW^{-1} \mY (\gamma - \Delta \vy),
\end{equation*}

the \textit{reduced KKT system} can be formulated as

\begin{equation}
    \begin{bmatrix}
        -\mH(\vx, \vy) & \mA(\vx)^\top \\
        \mA(\vx) & \mW \mY^{-1}
    \end{bmatrix}
    \begin{bmatrix}
        \Delta \vx \\
        \Delta \vy
    \end{bmatrix}
    = 
    \begin{bmatrix}
        \sigma \\
        \rho + \mW \mY^{-1} \gamma
    \end{bmatrix}.
\end{equation}

In practice, most interior algorithms solves the reduced KKT system iteratively for $\Delta \vx$ and $\Delta \vy$ and then updates the primal variables $\vx$, $\vw$, $\vy$ by a steplength $\alpha$:

\begin{equation}
    \begin{aligned}
        \vx^{(k+1)}&=\vx^{(k)}+\alpha^{(k)} \Delta \vx^{(k)}, \\
        \vw^{(k+1)}&=\vw^{(k)}+\alpha^{(k)} \Delta \vw^{(k)}, \\
        \vy^{(k+1)}&=\vy^{(k)}+\alpha^{(k)} \Delta \vy^{(k)}.
    \end{aligned}
\end{equation}

where the supscript $(k)$ indicates the $k$-th iteration. The steplength $\alpha^{(k)}$ is chosen according to certain criteria that is specific to each algorithm.

\subsection{A Simpler Case: Quadratic Programs}

Now consider a Quadratic Program (QP):

\begin{equation}
    \begin{aligned}
        \min &\quad \frac{1}{2}\vx^\top \mQ \vx + \vc^\top \vx \\
        \text{s.t.} &\quad \mA \vx \geq \vb\\
        &\quad \vx \geq \mathbf{0}.
    \end{aligned}
\end{equation}

To spedify the number of variables and constraints, it can be alternatively written as 

\begin{equation}
    \begin{aligned}
        \min & \quad \frac{1}{2}\vx^\top \mQ \vx + \vc^\top \vx \\
        \text{s.t.} & \quad \mA_j \vx \leq b_j,& \forall j \in \mathcal{M}\\
        &\quad x_i \geq 0, &\forall i \in \mathcal{N}.
    \end{aligned}
\end{equation}

where $\mathcal{M} = \{1, \ldots, m\}$ and $\mathcal{N} = \{1, \ldots, n\}$.

\subsubsection{Adding the Barrier Function}

According to the previous discussion, we can formulate the problem as

\begin{equation}
    \begin{aligned}
        \min &\quad \frac{1}{2}\vx^\top \mQ \vx + \vc^\top \vx - \mu \left[\ve^\top \log(\vw) + \ve^\top \log(\vx) \right].\\
        \text{s.t.} &\quad \mA \vx - \vw = \vb, \quad \vw \geq \mathbf{0}.
    \end{aligned}
\end{equation}

where $\log(\vx)$ is element-wise logarithm of $\vx$ and hence, $\ve^\top \log(\vx)$ is equivalent to $\sum_{i=1}^{n} \log(x_i)$. 

\subsubsection{The Lagrangian and KKT Conditions}

The Lagrangian for the problem is

\begin{equation}
    L(\vx, \vw, \vy) = \frac{1}{2}\vx^\top \mQ \vx + \vc^\top \vx - \mu \left[\ve^\top \log(\vw) + \ve^\top \log(\vx) \right] + \vy^\top(\vb -\mA \vx + \vw)
\end{equation}

By introducing $\vz = \mu\mX^{-1}\ve $ , the first-order KKT conditions can be written as

\begin{equation}
    \begin{aligned}
        \mA^\top \vy+\vz-\mQ \vx & = \vc \\
        \mA \vx-\vw & = \vb \\
        \mX \mZ \ve & =\mu \ve \\
        \mY \mW \ve & =\mu \ve .
    \end{aligned}
    \label{eq:QP-KKT}
\end{equation}

\subsubsection{Solving by Newton's Method}

In order to solve Equation \ref{eq:QP-KKT}, we can use the Newton's method by calculating $\Delta \vx$, $\Delta \vw$, $\Delta \vy$, $\Delta \vz$ for the following system:

\begin{equation}
    \begin{aligned}
        \mA^\top \Delta \vy+\Delta \vz-\mQ \Delta \vx & =\vc-\mA^\top \vy-\vz+\mQ \vx\\
        \mA \Delta \vx-\Delta \vw & =\vb-\mA \vx+\vw \quad\\
        \mZ \Delta \vx+\mX \Delta \vz+\Delta \mX \Delta \mZ \ve & =\mu \ve-\mX \mZ \ve \\
        \mW \Delta \vy+\mY \Delta \vw+\Delta \mY \Delta \mW \ve & =\mu \ve-\mY \mW \ve.
    \end{aligned}
\end{equation}

Notice that 

\begin{equation}
    \begin{aligned}
        \Delta \vz & =\mX^{-1}(\mu \ve-\mX \mZ \ve-\mZ \Delta \vx) \\
        \Delta \vw & =\mY^{-1}(\mu \ve-\mY \mW \ve-\mW \Delta \vy),
    \end{aligned}
\end{equation}

so what we actually need to solve is the following system:

\begin{equation}
    \begin{bmatrix}
    -\left(\mX^{-1} \mZ+\mQ\right) & \mA^\top \\
    \mA & \mY^{-1} \mW
    \end{bmatrix}
    \begin{bmatrix}
    \Delta \vx \\
    \Delta \vy
    \end{bmatrix}
    =
    \begin{bmatrix}
    \vc-\mA^\top \vy-\mu \mX^{-1} \ve+\mQ \vx \\
    \vb-\mA \vx+\mu \mY^{-1} \ve
    \end{bmatrix}.
    \label{eq:QP-reduced-KKT}
\end{equation}

Fortunately, in the case of quadratic programming, Equation \ref{eq:QP-reduced-KKT} has a closed-form solution

\begin{align}
    \Delta \vx &= \left(\mX^{-1} \mZ + \mQ \right)^{-1} \left[\mA^\top \Delta \vy - \vc - \mQ \vx + \mA^\top \vy + \mu \mX^{-1} \ve\right]\\
    \mP \Delta \vy &= \vb - \mA \vx + \mu \mY^{-1} \ve + \mA(\mX^{-1} \mZ + \mQ)^{-1} \left(\vc - \mA^\top \vy + \mQ \vx -\mu \mX^{-1} \ve\right)
\end{align}

where $\mP = \mA (\mX^{-1} \mZ + \mQ)^{-1} \mA^\top + \mY^{-1} \mW$. However, this expression is inefficient in practice since $(\mX^{-1} \mZ + \mQ)$ is not diagonal and hence difficult to calculate its inverse. Fortunately, we have an alternative

\begin{align}
    \Delta \vy &= \mY \mW^{-1} \left(\vb -\mA \vx + \mu \mY^{-1}\ve - \mA \Delta \vx \right)\\
    \mT \Delta \vx &= \vc - \mQ^\top\vy + \mQ \vx - \mu \mX^{-1} \ve - \mA^\top\mY\mW^{-1}(\vb - \mA\vx + \mu \mY^{-1}\ve)
\end{align}

where $\mT = -\left(\mX^{-1}\mZ + \mQ + \mA^\top \mY \mW^{-1}\mA\right)$. Notice that $\mT$ here is much easier to calculate than $\mP$. Then we calculate the steplength $\alpha$ by some rule and do the following update:

\begin{equation}
    (\vx, \vw, \vy, \vz) \leftarrow (\vx, \vw, \vy, \vz) + \alpha \left(\Delta \vx, \Delta \vy, \Delta \vz, \Delta \vw\right)
\end{equation}

where the supscript $(k)$ is omitted for simplicity. For the steplength $\alpha$, there are different rules to choose it. A rule of the thumb is to find the largest $\alpha$ such that

\begin{equation}
    \min_{i,j} \left\{ (\vx +\alpha \Delta \vx)_i (\vz +\alpha \Delta \vz)_i, (\vy +\alpha \Delta \vy)_j (\vw +\alpha \Delta \vw)_j \right\} \geq 0
    \label{eq:practical-steplength-rule}
\end{equation}

and then

\begin{equation}
    (\vx, \vw, \vy, \vz) \leftarrow (\vx, \vw, \vy, \vz) + 0.99 \alpha \left(\Delta \vx, \Delta \vy, \Delta \vz, \Delta \vw\right)
\end{equation}

\subsubsection{A Practical Algorithm}

Now we present a formal algorithm used in practice \citep{vanderbeiLOQO1999}. This algorithm effectively combines the conjugate gradient method and the IPM for solving QPs.

\begin{figure}[H]
    \begin{minipage}[t]{0.48\textwidth}
        \begin{algorithm}[H]
            \small
            \caption{Conjugate Gradient Method for Equation \ref{eq:QP-reduced-KKT}}
            \label{alg:CG}
            \begin{algorithmic}[1]
                \REQUIRE 
                \STATE $\vp\leftarrow \vb - \mA \vx + \mu \mY^{-1} \ve + \mA(\mX^{-1} \mZ + \mQ)^{-1}$ \\
                $\left(\vc - \mA^\top \vy + \mQ \vx -\mu \mX^{-1} \ve\right)$ \label{line:CG-1}
                \STATE $\vv \leftarrow -\vp$ \label{line:CG-2}
                \STATE $\Delta \vy \leftarrow \mathbf{0}$ \label{line:CG-3}
                \FOR{$m$ iterations}
                    \STATE $\vu \leftarrow \mA (\mX^{-1} \mZ + \mQ)^{-1} \mA^\top \vp + \mY^{-1} \mW \vp$ \label{line:CG-5}
                    \STATE $\alpha \leftarrow \vv^\top \vv / \vp^\top \vu$ \label{line:CG-6}
                    \STATE $\Delta \vy \leftarrow \Delta \vy + \alpha \vp$ \label{line:CG-7}
                    \STATE $\vv_{\text{new}} \leftarrow \vv + \alpha \vu$ \label{line:CG-8}
                    \STATE $\beta \leftarrow \vv_{\text{new}}^\top \vv_{\text{new}} / \vv^\top \vv$ \label{line:CG-9}
                    \STATE $\vv \leftarrow \vv_{\text{new}}$ \label{line:CG-10}
                    \STATE $\vp \leftarrow -\vv + \beta \vp$ \label{line:CG-11}
                \ENDFOR
                \RETURN $\Delta \vy$ that solves Equation \ref{eq:QP-reduced-KKT}
            \end{algorithmic}
        \end{algorithm}
    \end{minipage}
    \hfill
    \begin{minipage}[t]{0.48\textwidth}
        \begin{algorithm}[H]
            \small
            \caption{A Pratical Interior-Point Algorithm for Quadratic Programming}\label{alg:IPM}
            \begin{algorithmic}[1]
                \REQUIRE Initial values of $\vx$, $\vw$, $\vy$, $\vz$, hyperparameter $\delta \in (0,1)$
                \STATE $\mu = \left(\vz^\top \vx + \vy^\top\vw\right) / (n+m)$ \label{line:IPM-1}
                \REPEAT
                    \STATE Solve Equation \ref{eq:QP-reduced-KKT} to get $\Delta \vy$ \label{line:IPM-3}
                    \STATE $\Delta \vx \leftarrow \left(\mX^{-1} \mZ + \mQ\right)^{-1}$ \\
                    $\left[\mathbf{A}^\top \Delta \mathbf{y} - \mathbf{c} - \mathbf{Q} \mathbf{x} + \mathbf{A}^\top \mathbf{y} + \mu \mathbf{X}^{-1} \mathbf{e}\right]$ \label{line:IPM-4}
                    \STATE $\Delta \vz \leftarrow \mX^{-1} (\mu \ve - \mX \mZ \ve - \mZ \Delta \vx)$ \label{line:IPM-5}
                    \STATE $\Delta \vw \leftarrow \mY^{-1} (\mu \ve - \mY \mW \ve - \mW \Delta \vy)$ \label{line:IPM-6}
                    \STATE Calculate $\alpha$ according to Equation \ref{eq:practical-steplength-rule} \label{line:IPM-7}
                    \STATE Update $(\vx, \vw, \vy, \vz) \leftarrow (\vx, \vw, \vy, \vz) + 0.99 \alpha \left(\Delta \vx, \Delta \vy, \Delta \vz, \Delta \vw\right)$ \label{line:IPM-8}
                    \STATE $\mu \leftarrow \delta \mu$ \label{line:IPM-9}
                \UNTIL{optimality condition is satisfied}
            \end{algorithmic}
            \vspace{0.6em}
        \end{algorithm}
    \end{minipage}
\end{figure}

In order to solve Equation \ref{eq:QP-reduced-KKT}, we can use a \textit{conjugate gradient} method as is detailed in Algorithm \ref{alg:CG}. The conjugate gradient method is an iterative algorithm for solving linear systems that is particularly useful for large, sparse systems arising in various optimization problems. 

In Algorithm \ref{alg:CG}, we initialize the search direction \(\vp\) using the given values of \(\vx\), \(\vy\), \(\vz\), and \(\vw\). The variable \(\vv\) is then set to the negative of \(\vp\), and \(\Delta \vy\) is initialized to zero. In each iteration, we compute the search direction \(\vu\), the step size \(\alpha\), and update \(\Delta \vy\), \(\vv\), and \(\vp\). The process repeats for a specified number of iterations or until convergence criteria are met, yielding \(\Delta \vy\) that solves Equation \ref{eq:QP-reduced-KKT}.

Algorithm \ref{alg:IPM} describes a practical IPM for quadratic programming. This method starts with initial values for the variables \(\vx\), \(\vw\), \(\vy\), and \(\vz\), and a hyperparameter \(\delta\). At each iteration, the algorithm solves the reduced KKT system to obtain \(\Delta \vy\), updates the primal and dual variables using \(\Delta \vx\), \(\Delta \vy\), \(\Delta \vz\), and \(\Delta \vw\), and adjusts the barrier parameter \(\mu\). The steplength \(\alpha\) is calculated to ensure feasibility and convergence of the solution.

The combination of these algorithms ensures robust and efficient solving of QPs, demonstrating the practical applicability of IPMs.

\subsection{Proof of the Equivalence of MPNNs and IPMs}
\label{appendix:proof_equivalence_of_MPNNs_and_IPMs}

In this section, we will show that MPNNs are equivalent to the general framework of IPMs for solving quadratic programming problems. We then show that our UniEGNN fits into the framework of MPNNs. To do this, we first show that there exists a MPNN $f_{\text{CG}}$ that reproduces the conjugate gradient method in Algorithm \ref{alg:CG} (Lemma \ref{lemma:CG}). We then show that there exists a MPNN $f_{\text{IPM}}$ that reproduces the IPM in Algorithm \ref{alg:IPM} (Theorem \ref{theorem:UniEGNN_IPM}). Finally, UniEGNN belongs to the class of MPNNs, which completes the proof.

\textbf{Proof of Lemma \ref{lemma:CG}}. There exists a MPNN \( f_{\text{CG}} \) composed of a finite number of message-passing steps that reproduces the conjugate gradient method described in Algorithm \ref{alg:CG}. Specifically, for any QP instance \( I = (\mQ, \vc, \mA, \vb) \) and any initial point \((\vx, \vy, \vz, \vw) > 0\), \( f_{\text{CG}} \) maps the graph \( \mathcal{H}(I) \) carrying the initial values on the variable and constraint nodes to the graph \( \mathcal{H}(I) \) carrying the output \(\Delta \vy\) of Algorithm \ref{alg:CG} on the constraint nodes.

\begin{proof}
    We prove by showing that each line in Algorithm \ref{alg:CG} can be reproduced by several message-passing steps within the UniEGNN framework.
    
    \begin{itemize}[left=0pt, topsep=0pt, itemsep=0pt]
        \item For Line \ref{line:CG-1}, the computation can be broken down as follows. First, we compute \(\vh_1 = \mA^\top \vy\) by a constraints-to-variables message-passing step. Next, we compute \(\vh_2 = \mu \ve\) by an objective-to-variables message-passing step, followed by \(\vh_3 = \vc\) and \(\vh_4 = \mQ \vx\) also by objective-to-variables message-passing steps. Then, \(\vh_5 = \mQ\) is computed similarly. We then perform a local operation on variable nodes to compute \(\vh_6 = -\vx + (\mX^{-1}\mZ + \vh_5)(\vh_3 - \vh_1 + \vh_4 - \vh_2 \mX)\). Subsequently, \(\vh_7 = \mA \vh_6\) is obtained by a variables-to-constraints message-passing step, followed by \(\vh_8 = \mu \ve\) computed by an objective-to-constraints message-passing step. Finally, \(\vp = \vb + \vh_7 - \vh_8 \mY^{-1}\) is computed as a local operation on constraint nodes.

        \item Line \ref{line:CG-2} is a local operation on constraint nodes: \(\vv = -\vp\).

        \item Line \ref{line:CG-3} is a local operation on constraint nodes: \(\Delta \vy = 0\).

        \item For Line \ref{line:CG-5}, the computation is broken down as follows. We first compute \(\vh_1 = \mA^\top \vp\) by a constraints-to-variables message-passing step, followed by \(\vh_2 = \mQ\) by an objective-to-variables message-passing step. Next, we compute \(\vh_3 = (\mX^{-1}\mZ + \vh_2)\vh_1\) as a local operation on variable nodes, then \(\vh_4 = \mA \vh_3\) by a variables-to-constraints message-passing step. Finally, \(\vu = \vh_4 + \mY^{-1}\mW \vp\) is computed as a local operation on constraint nodes.

        \item For Line \ref{line:CG-6}, the computation is broken down as follows. Compute \(\vh_1 = \vv^\top \vv\) by a constraints-to-objective message-passing step, followed by \(\vh_2 = \vp^\top \vu\) by another constraints-to-objective message-passing step. Finally, compute \(\alpha = \vh_1 / \vh_2\) as a local operation on the objective node.

        \item For Line \ref{line:CG-7}, the computation is broken down as follows. First, compute \(\vh_1 = \alpha \vp\) by an objective-to-constraints message-passing step, then update \(\Delta \vy = \Delta \vy + \vh_1\) as a local operation on constraint nodes.

        \item For Line \ref{line:CG-8}, the computation is broken down as follows. First, compute \(\vh_1 = \alpha \vu\) by an objective-to-constraints message-passing step, then update \(\vv_{\text{new}} = \vv + \vh_1\) as a local operation on constraint nodes.

        \item For Line \ref{line:CG-9}, the computation is broken down as follows. First, compute \(\vh_1 = \vv_{\text{new}}^\top \vv_{\text{new}}\) by a constraints-to-objective message-passing step, followed by \(\vh_2 = \vv^\top \vv\) by another constraints-to-objective message-passing step. Finally, compute \(\beta = \vh_1 / \vh_2\) as a local operation on the objective node.

        \item Line \ref{line:CG-10} is a local operation on constraint nodes: update \(\vv = \vv_{\text{new}}\).

        \item For Line \ref{line:CG-11}, the computation is broken down as follows. First, compute \(\vh_1 = \beta \vp\) by an objective-to-constraints message-passing step, then update \(\vp = -\vv + \vh_1\) as a local operation on constraint nodes.
    \end{itemize}
\end{proof}

\textbf{Proof of Theorem \ref{theorem:UniEGNN_IPM}}. There exists an MPNN \( f_{\text{IPM}} \) composed of a finite number of message-passing steps that reproduces each iteration of the IPM described in Algorithm \ref{alg:IPM}. Specifically, for any QP instance \( I = (\mQ, \vc, \mA, \vb) \) and any iteration step \( t \geq 0 \), \( f_{\text{IPM}} \) maps the hypergraph \( \mathcal{H}(I) \) carrying the current iterate values \((\vx_{t}, \vz_{t}), (\vy_{t},\vw_{t})\) on the variable and constraint nodes to the hypergraph \( \mathcal{H}(I) \) carrying the next iterate values \((\vx_{t+1}, \vz_{t+1}), (\vy_{t+1},\vw_{t+1})\) on the variable and constraint nodes.

\begin{proof}
    Each step in Algorithm \ref{alg:IPM} can be reproduced by several message-passing steps within the UniEGNN framework.

    \begin{itemize}[left=0pt, topsep=0pt, itemsep=0pt]
        \item For Line \ref{line:IPM-1}, the computation is broken down as follows. Compute \(\vh_1 = \vz^\top \vx\) by a variable-to-objective message-passing step, followed by \(\vh_2 = \vy^\top \vw\) by a constraints-to-objective message-passing step. Finally, compute \(\mu = (\vh_1 + \vh_2) / (n+m)\) as a local operation on the objective node.

        \item For Line \ref{line:IPM-3}, the computation is broken down as follows. First, compute \(\vh_1 \leftarrow \mA^\top (\vy + \Delta \vy)\) by a constraints-to-variables message-passing step. Then, compute \(\vh_2 \leftarrow \mu \ve\), \(\vh_3 \leftarrow \vc\), and \(\vh_4 \leftarrow \mQ \vx\) by objective-to-variables message-passing steps. Compute \(\vh_5 \leftarrow \mQ\) similarly. Perform a local operation on variable nodes to compute \(\vh_6 \leftarrow \vh_5 + \mX^{-1}\mZ\) and \(\vh_7 \leftarrow \vh_6^{-1}\). Finally, compute \(\Delta \vx \leftarrow \vh_7(\vh_1 + \vh_2 \mX - \vh_3 - \vh_4)\) as a local operation on variable nodes.

        \item For Line \ref{line:IPM-4}, the computation is broken down as follows. Compute \(\vh_1 \leftarrow \mu \ve\) by an objective-to-variables message-passing step, followed by \(\Delta \vz \leftarrow \mX^{-1}\vh_1 - \vz - \mX^{-1}\mZ\Delta \vx\) as a local operation on variable nodes.

        \item For Line \ref{line:IPM-5}, the computation is broken down as follows. Compute \(\vh_1 \leftarrow \mu \ve\) by an objective-to-constraints message-passing step, followed by \(\Delta \vw \leftarrow \mY^{-1} \vh_1 - \vw - \mY^{-1}\mW \Delta \vy\) as a local operation on constraint nodes.

        \item For Line \ref{line:IPM-6}, the computation can be performed by message-passing steps as follows. On every variable node, compute:
        \begin{equation*}
            \alpha_i = \max \left\{\alpha \in(0, \infty) \mid \alpha^2 \Delta x_i \Delta z_i + \alpha \left(x_i \Delta z_i + \Delta x_i z_i \right) + x_i z_i \geq 0 \right\}
        \end{equation*}
        On every constraint node, compute:
        \begin{equation*}
            \tilde{\alpha}_j = \max \left\{\alpha \in(0, \infty) \mid \alpha^2 \Delta y_i \Delta w_i + \alpha \left(y_i \Delta w_i + \Delta y_i w_i \right) + y_i w_i \geq 0 \right\}
        \end{equation*}
        Then, the computation of $\alpha_v$ and $\tilde{\alpha}_c$ can be viewed as variable-to-objective and constraint-to-objective message-passing steps, respectively:
        \begin{equation*}
            \alpha_v \leftarrow \min_i \alpha_i, \quad \tilde{\alpha}_c \leftarrow \min_j \tilde{\alpha}_j
        \end{equation*}
        Finally, $\alpha$ can be computed as local operations on the objective node:
        \begin{equation*}
            \alpha \leftarrow \min \left(\alpha_v, \tilde{\alpha}_c\right)
        \end{equation*}

        \item For Line \ref{line:IPM-7}, the computation is broken down as follows. Compute \(\vh_1 \leftarrow \alpha \ve\) by an objective-to-variables message-passing step. Update \(\vx \leftarrow \vx + 0.99 \mD(\vh_1) \Delta \vx\) and \(\vz \leftarrow \vz + 0.99 \mD(\vh_1) \Delta \vz\) as local operations on variable nodes. Compute \(\vh_2 \leftarrow \alpha \ve\) by an objective-to-constraints message-passing step. Update \(\vy \leftarrow \vy + 0.99 \mD(\vh_2) \Delta \vy\) and \(\vw \leftarrow \vw + 0.99 \mD(\vh_2) \Delta \vw\) as local operations on constraint nodes.

        \item Line \ref{line:IPM-8} is a local operation on the objective node: update \(\mu \leftarrow \delta \mu\).
    \end{itemize}
\end{proof}

\subsection{Conclusion}

In this section, we have demonstrated that the UniEGNN framework can reproduce the steps of the IPM for solving QPs. This establishes a theoretical foundation for the equivalence between UniEGNN and IPMs, thereby validating the applicability of UniEGNN. While it is challenging to explicitly formulate the IPM for general nonlinear optimization, the UniEGNN framework can potentially be extended to handle these more complicated cases.

\section{Benchmark Problems and Datasets}
\label{appendix:benchmark_problems}

In Section \ref{sec:experiments}, we evaluate our framework NeuralQP on two benchmark problems: the Quadratic Multiple Knapsack Problem (QMKP) and the Random Quadratically Constrained Program (RandQCP). These benchmarks are chosen due to their relevance in testing the efficiency and effectiveness of optimization algorithms under complex, real-world constraints. In the following subsections, we provide a detailed description of each benchmark problem, including their formulation, characteristics, and the datasets used for experimentation.

\subsection{QMKP}
\label{appendix:QMKP}

The Quadratic Multiple Knapsack Problem (QMKP) \citep{hiley2006quadratic} is a combination of the Quadratic Knapsack Problem (QKP) \citep{gallo1980quadratic} and the Multiple Knapsack Problem (MKP) \citealp{kellerer2004multidimensional}, adapted to better reflect real-world scenarios. QKP assumes a single-dimensional restriction on the items, which has limited applications; MKP assumes a linear relationship between the values of items, which ignores the interactions between them. To overcome these limitations, the QMKP incorporates not only a linear value for each item but also quadratic terms representing the interactions between items when selected together, as well as multidimensional constraints. This additional complexity allows for more sophisticated modeling and is applicable in various fields. The QMKP can be mathematically formulated as follows:

\begin{equation}
\begin{aligned}
\max\quad &\sum_{i} c_i x_i + \sum_{(i,j)\in E} q_{ij}x_i x_j,\\
\text{s.t.}\quad &a_i^k x_i \leq b^k, \quad \forall k \in M,\\
&x_i \in \mathbb{B}, \quad \forall i \in N,
\end{aligned}
\end{equation}

where \(c_i\) represents the value of item \(i\), and the coefficient \(q_{ij}\) denotes the interaction between items \(i\) and \(j\) when selected together. The coefficients \(a_i^k\) relate to the attribute of item \(i\) concerning constraint \(k\), and \(b^k\) is the upper limit for that attribute. In our experiment, $c_i$, $q_{ij}$ and $a_i^k$ are generated from the uniform distribution $U(0,1)$, with $b^k = \frac{1}{2}\sum_{i\in N} a_i^k$ so that the 0s and 1s in the optimal solution are roughly balanced and that the constraints remain valid.

The QMKP thus extends the scope of the QKP and the MKP, making it more suitable for scenarios where the value of items depends on the simultaneous selection of others.  For example, consider optimizing a day's diet to maximize the nutritional value of each food item. In this scenario, certain food combinations may offer greater nutritional benefits, akin to the standard QKP. However, each food item also possesses various attributes such as calorie, carbohydrate, and sodium content. To ensure a healthy diet, these attributes must not exceed specific limits. The QMKP's ability to capture the interdependencies among items makes it valuable for addressing complex optimization challenges in real-world settings with intricate constraints and objectives.

\subsection{RandQCP}
\label{appendix:RandQCP}

The Random Quadratically Constrained Quadratic Program (RandQCP) is an extension of the Independent Set (IS) problem \citep{coxeter1950self}, where we introduce quadratic non-convex constraints to enhance its difficulty. In a graph, we represent distinct vertex attributes by introducing weights \(c_i\) into the objective function, where each weight corresponds to the attribute of the node \(i\). To capture high-order relationships, we employ a hypergraph model, denoted as \(\mathcal{H}=(\mathcal{V}, \mathcal{E})\), where \(\mathcal{V}\) is the set of vertices and \(\mathcal{E}\) is the set of hyperedges. In such a hypergraph model, each hyperedge can encompass multiple vertices, reflecting the multi-dimensional interactions within each hyperedge. The primary objective of RandQCP is to identify the largest weighted set of vertices that satisfies the constraint within each hyperedge. Therefore, the mathematical form of RandQCP can be expressed as:

\begin{equation}
\begin{aligned}
    \max\quad &\sum_{i\in V}c_i x_i,\\
    \text{s.t.}\quad &\sum_{i \in e}a_i x_i + \sum_{i,j \in e, i\neq j} q_{ij}x_i x_j - |e|\leq 0,\quad \forall e \in \mathcal{E},\\
    &x_i \in \mathbb{B}, \quad \forall i \in V.
\end{aligned}
\end{equation}

The coefficients \(a_i\) and \(q_{ij}\) are randomly generated from a uniform distribution \(U(0,1)\). The term \(|e|\) represents the degree (i.e., the number of vertices) of the hyperedge, which is chosen to ensure that the constraints are both valid and feasible.

\subsection{Dataset Generation}
\label{appendix:dataset_generation}

Since existing datasets and learning methods couldn't meet the training and testing requirements (Appendix \ref{appendix:solution_strategy}), we generated the two benchmark problems using random coefficients described in Appendix \ref{appendix:QMKP} and \ref{appendix:RandQCP}. The numbers of variables and constraints for the problems are shown in Table \ref{table:benchmark_problems}. For training datasets, we generated 10800 instances for Mini-scale problems, 1000 instances for 1000-scale problems; 100 instances for 2000-scale problems. We used Gurobi to solve the problems until the gap was less than or equal to 10\%, and we considered the solution at this point as the optimal or near-optimal solution for the training set, since the rest of the time would be primarily spent on improving dual bounds. For testing datasets, we generated 100 instances for 1000-scale problems, and 10 instances for 2000-scale, 5000-scale, and 10000-scale problems, respectively. We did not generate 5000-scale and 10000-scale instances for training since solving these large-scale instances to optimality is computationally expensive, and such computation burden is precisely what we aim to overcome with our framework. Additionally, we did not generate Mini-scale instances for testing since our framework is targeted at large-scale problems.

\begin{table}[h]
    \centering
    \small
    \renewcommand{\arraystretch}{1.2}
    \caption{Details of the Two Benchmark Problems. ``-'' means no instances are generated.}
    \label{table:benchmark_problems}
    \begin{tabular}{cccccc}
    \hline
    \multirow{2}{*}{\textbf{Problem}} &
      \multirow{2}{*}{\textbf{Scale}} &
      \multirow{2}{*}{\textbf{\begin{tabular}[c]{@{}c@{}}Number of\\ Variables\end{tabular}}} &
      \multirow{2}{*}{\textbf{\begin{tabular}[c]{@{}c@{}}Number of\\ Constraints\end{tabular}}} &
      \multicolumn{2}{c}{\textbf{Number of Instances}} \\ \cline{5-6} 
     &       &                      &      & \textbf{Train} & \textbf{Test} \\ \hline
    \multirow{7}{*}{\begin{tabular}[c]{@{}c@{}}QMKP\\ (Max)\end{tabular}} &
      \multirow{3}{*}{Mini} &
      400 &
      5 &
      3600 &
      - \\
     &       & 500                  & 5    & 3600           & -             \\
     &       & 600                  & 5    & 3600           & -             \\ \cline{2-6} 
     & 1000  & 1000                 & 5    & 1000           & 100           \\
     & 2000  & 2000                 & 10   & 100            & 10            \\
     & 5000  & 5000                 & 10   & -              & 10            \\
     & 10000 & 10000                & 20   & -              & 10            \\ \hline
    \multirow{10}{*}{\begin{tabular}[c]{@{}c@{}}RandQCP\\ (Max)\end{tabular}} &
      \multirow{6}{*}{Mini} &
      \multirow{2}{*}{100} &
      60 &
      1800 &
      - \\
     &       &                      & 80   & 1800           & -             \\
     &       & \multirow{2}{*}{200} & 120  & 1800           & -             \\
     &       &                      & 160  & 1800           & -             \\
     &       & \multirow{2}{*}{300} & 180  & 1800           & -             \\
     &       &                      & 240  & 1800           & -             \\ \cline{2-6} 
     & 1000  & 1000                 & 800  & 1000           & 100           \\
     & 2000  & 2000                 & 1800 & 100            & 10            \\
     & 5000  & 5000                 & 4000 & -              & 10            \\
     & 10000 & 10000                & 8000 & -              & 10            \\ \hline
    \end{tabular}
    \end{table}
\section{Main Experiment Details}
\label{appendix:experiments_details}

\subsection{Baselines}
\label{appendix:baselines}
We use the currently most advanced and widely used solvers, Gurobi and SCIP, as baselines. The Gurobi version used is 11.0.1 with \verb|Threads| set to 4, \verb|NonConvex| set to 2, and other parameters set to default values. The SCIP version used is 4.3.0 with all default settings. The scaled-constrained versions are restricted to 30\% and 50\% of the number of variables in the original problem. To ensure the validity of comparative experiments, the scale-constrained versions share the same parameters as the no-constrained versions.

\subsection{UniEGNN Model Formulation}
\label{appendix:model_formulation}

Given a QCQP instance $I$ of the form in Equation \ref{eq:qp} with binary variables, the optimal solution $\mathcal{X}^*$ can be viewed as a function $g$ of the instance $I$, denoted as $g(I)$. The learning objective is then to approximate this function $g$ using a neural network model $g_\theta$, parameterized by $\theta$. Formally, the training dataset is constructed as $\mathcal{D}=\{(I_i,\mathcal{X}_i^*)\}_{i=1}^N$, where $\{I_i\}_{i=1}^N$ represents $N$ problem instances and $\{\mathcal{X}_i^*\}_{i=1}^N$ denotes the corresponding optimal or near-optimal solutions obtained by off-the-shelf solvers.

To measure the difference between the network's predictions and the true solutions, we employ the binary cross-entropy loss with logits, represented by \verb|BCEWithLogitLoss| in \verb|PyTorch|. This loss function is suitable because the network outputs are not probabilities. A sigmoid function is applied to convert them into probabilities before calculating the loss. The loss for each instance in the dataset is calculated as:
\begin{equation}
    \mathcal{L}(\theta) = - \frac{1}{N} \sum_{i=1}^N \left[ y_i \cdot \log(\sigma(g_\theta(I_i))) + (1 - y_i) \cdot \log(1 - \sigma(g_\theta(I_i))) \right],
    \label{eq:neural_network}
\end{equation}
where $\sigma$ denotes the sigmoid function, $y_i$ is the true label indicating whether the solution is 0 or 1, and $g_\theta(I_i)$ is the raw output of the network for the given QCQP instance $I_i$. The goal during training is to minimize this loss function, effectively adjusting the parameters $\theta$ to improve the approximation of the true function $g$.

\subsection{General Experiment settings}
\label{appendix:experiment_settings}

In the training stage, we implemented a UniEGNN model with 6 convolution layers and an MLP model with 3 layers. For each combination of problem type and scale, a UniEGNN model was trained on a single GPU of a machine with 4 NVIDIA Tesla V100(32G) GPUs for 100 epochs. The dimensions of the initial embedding space, the hidden space, the MLP input (also the convolution output), and the final output were 16, 64, 16, and 1 respectively. The number of convolutional layers of different models was all set to 6. The first and second aggregation functions for the convolutional layers were SUM and MEAN respectively (see Section \ref{sec:neural_prediction}). The optimizer we adopted was \verb+AdamW+ from \verb+torch.optim+ with \verb+lr=1e-4+ and \verb+weight_decay=1e-4+. The activation we used was \verb+leakyrelu+ with \verb+negative_slope=0.1+.

In the testing stage, experiments were conducted on a machine with 128 Intel(R) Xeon(R) Platinum 8375C CPU @ 2.90GHz CPUs. During the initial feasible solution search phase, the small-scale solver was configured to return when one feasible solution is found or time limit is reached. In the Parallel Neighborhood Optimization phase, the small-scale solver was set to return when the optimal solution or the time limit is reached. In experiments comparing runtime, we set the objective value limits to be the minimum value among the 10 instances obtained using NeuralQP in experiments comparing objective value. Each value in Table \ref{table:2} is averaged among 10 instances. If some instances didn't reach the objective value given before the time limit is reached, we use the time limit given as the result to calculate the average and add $>$ before the values. Besides, we limited the initial feasible solution search, each round of neighborhood search, and neighborhood crossover time to 30s for small-scale, 2000-scale, and 5000-scale problems and 60s for 10000-scale problems.

\subsection{Convergence Analysis}
\label{appendix:convergence_analysis}
\begin{figure}[h]
  \begin{minipage}[H]{0.50\textwidth}
    \centering
    \includegraphics[width=\textwidth]{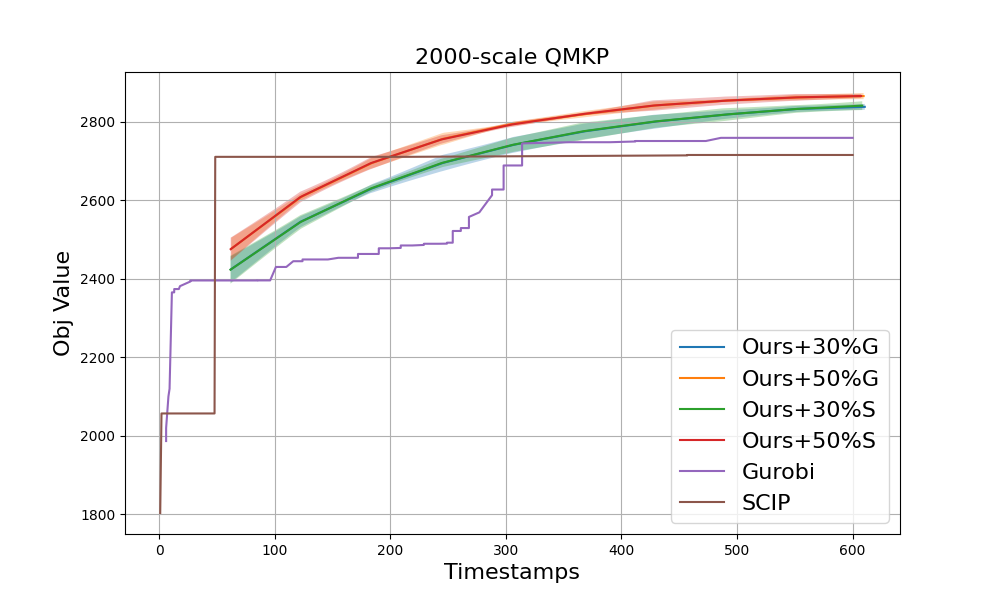}
  \end{minipage}
  \begin{minipage}[H]{0.50\textwidth}
    \centering
    \includegraphics[width=\textwidth]{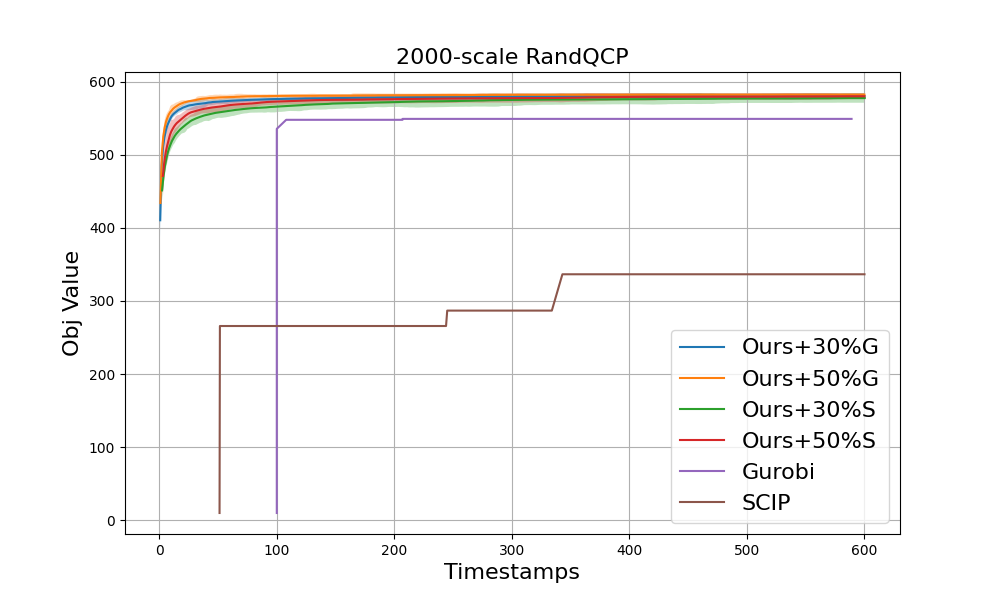}
  \end{minipage}
  \caption{The time-objective value figure for 2000-scale problems.}
  \label{fig:2000-scale}

  \begin{minipage}[H]{0.50\textwidth}
    \centering
    \includegraphics[width=\textwidth]{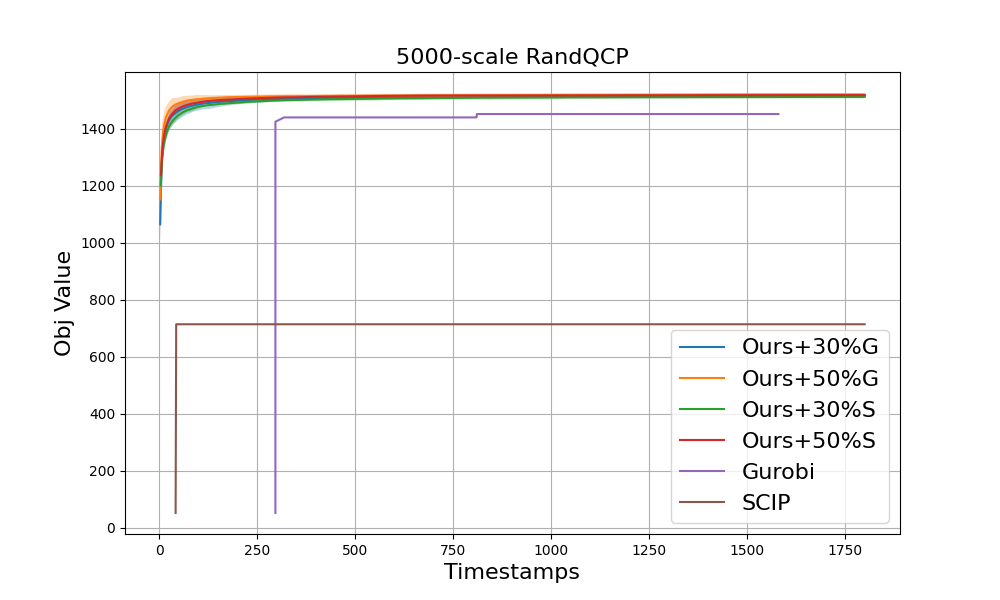}
  \end{minipage}
  \begin{minipage}[H]{0.50\textwidth}
    \centering
    \includegraphics[width=\textwidth]{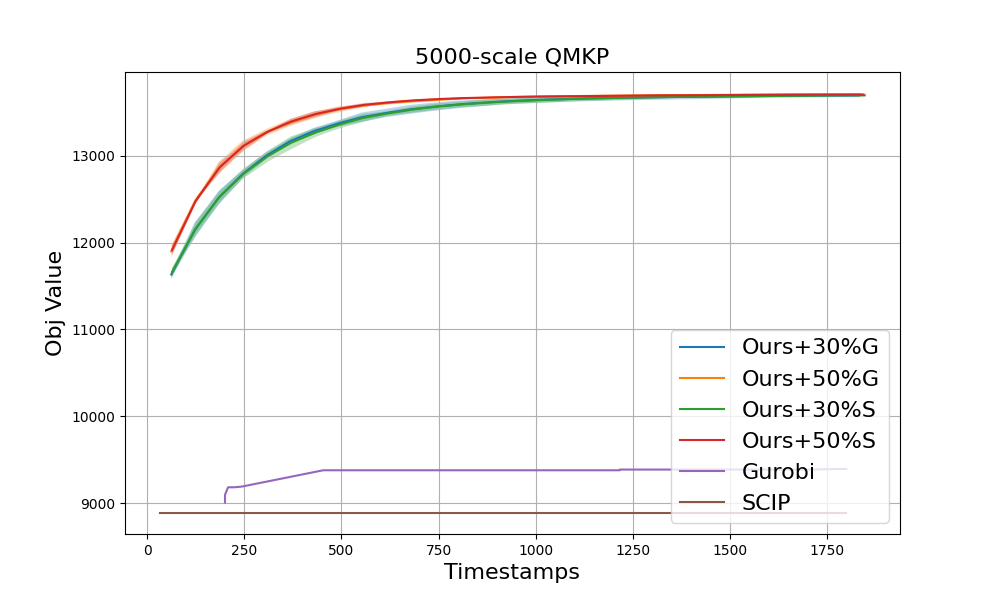}
  \end{minipage}
  \caption{The time-objective value figure for 5000-scale problems.}
  \label{fig:5000-scale}

    \begin{minipage}[H]{0.50\textwidth}
    \centering
    \includegraphics[width=\textwidth]{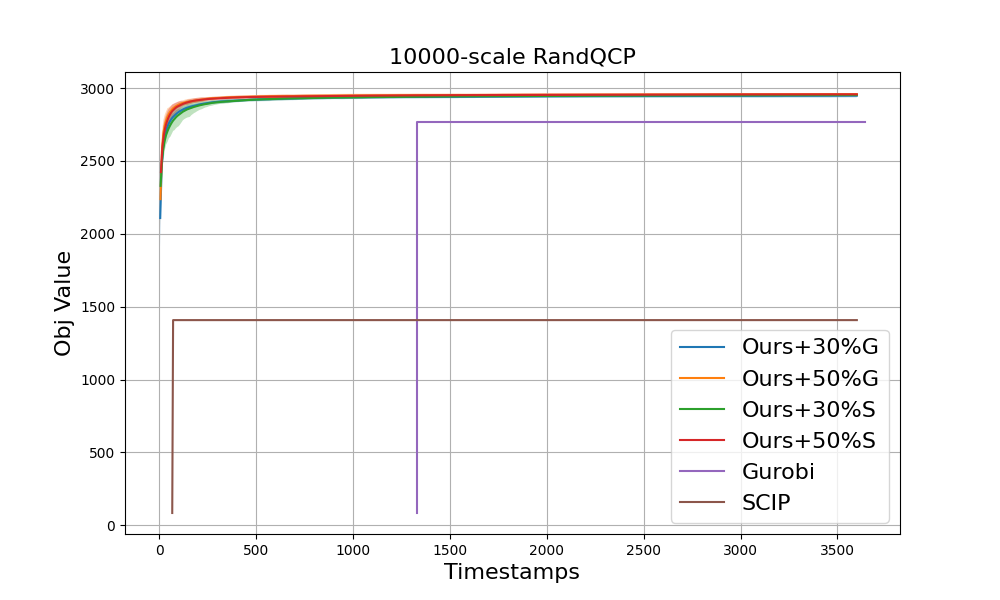}
  \end{minipage}
  \begin{minipage}[H]{0.50\textwidth}
    \centering
    \includegraphics[width=\textwidth]{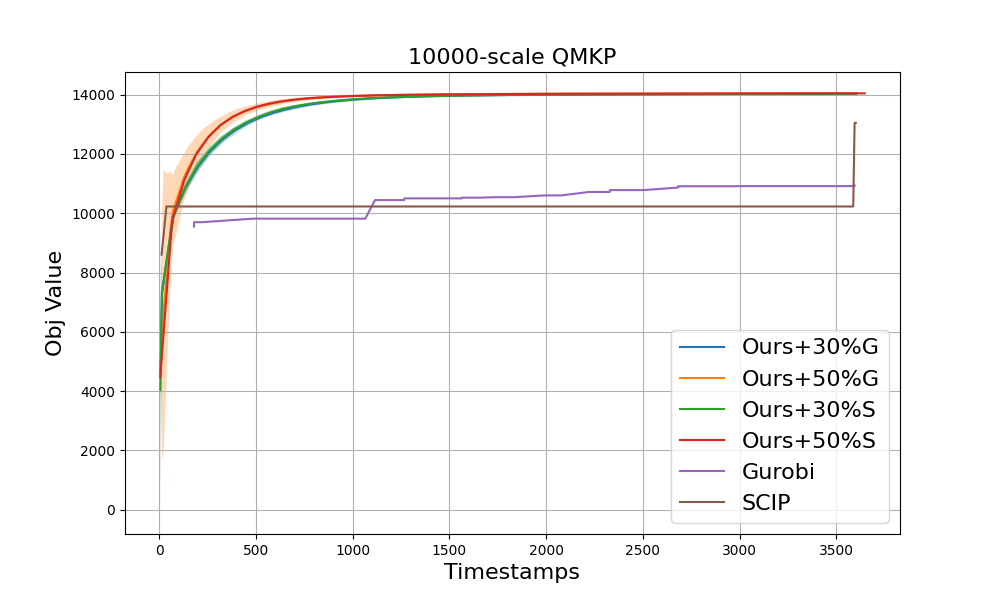}
  \end{minipage}
  \caption{The time-objective value figure for 10000-scale problems.}
  \label{fig:10000-scale}
\end{figure}
Convergence is an important criterion for assessing the effectiveness of neighborhood optimization. It can be obtained by observing the curve of the objective value over time during the optimization process. The convergence curves for solving QCQPs by using our framework with the small-scale version solver, in comparison with the large-scale version solver, are shown in Figure \ref{fig:2000-scale}, Figure \ref{fig:5000-scale} and Figure \ref{fig:10000-scale}. We ran our method on each instance for 5 times. The solid line represents the mean values, while the shaded area represents the $3\sigma$ error bar (mean value plus/minus the standard deviation multiplied by 3).

Our framework exhibits superior convergence performance, especially for large-scale problems, achieving high-quality solutions in less time compared to Gurobi and SCIP. Additionally, the shaded area indicates consistent results across multiple runs, demonstrating NeuralQP's superior stability.

\section{Extra Experiments}
\subsection{Ablation Studies}
\label{appendix:ablation_studies}

We further conducted ablation studies on each of the two stages of our framework, i.e., Parallel Neighborhood Optimization and Neural Prediction, to further validate the effectiveness and generalization ability of our framework.

\subsubsection{Neural Prediction}
\label{appendix:neural_prediction}
    
We tested our method in the Neural Prediction stage against the state-of-the-art solver Gurobi. We trained multiple UniEGNN models on instances of varying scales and then tested the model on instances of varying scales. To be specific, we let our model predict a solution and then use the initial feasible solution strategy to acquire a feasible solution. For Gurobi, we set \verb+SolutioniLimit=2+ since the first solution is trivial and has the objective value of 0 which is meaningless; we set\verb+Timelimit=10+ to prevent gurobi from running too long. We ran both our method and Gurobi on each instance for 1 time. Since our method achieves close objective values and has similar runtime among different runs on a single instance, we report the mean objective value and runtime on 10 instances of each method in Table \ref{tab:ablation_neural_prediction}. We do not report the std due to the limitation of space. Please refer to supplementary materials for complete raw data.

\begin{table}[H]
\renewcommand{\arraystretch}{1.3}
\caption{Comparison of Objective Value and Runtime for Neural Prediction}
\label{tab:ablation_neural_prediction}
\resizebox{\textwidth}{!}{%
\begin{tabular}{ccccccccccc}
\hline
\multirow{2}{*}{Problem} &
  \multirow{2}{*}{Train} &
  \multirow{2}{*}{Test} &
  \multicolumn{2}{c}{Objective Value} &
  \multicolumn{2}{c}{Runtime} &
  \multirow{2}{*}{Obj Ratio} &
  \multirow{2}{*}{Runtime Ratio} &
  \multirow{2}{*}{Faster} &
  \multirow{2}{*}{Better} \\ \cline{4-7}
                          &                       &       & Ours     & Gurobi  & Ours & Gurobi &       &      &     &     \\ \hline
\multirow{12}{*}{RandQCP} & \multirow{4}{*}{1000} & 1000  & 185.72   & 9.24    & 0.11 & 0.02   & 20.10 & 5.80 & 0   & 100 \\
                          &                       & 2000  & 310.61   & 10.89   & 0.24 & 0.08   & 28.54 & 3.11 & 0   & 10  \\
                          &                       & 5000  & 931.28   & 44.56   & 0.55 & 0.18   & 20.90 & 3.11 & 0   & 10  \\
                          &                       & 10000 & 1860.34  & 93.57   & 1.08 & 0.34   & 19.88 & 3.15 & 0   & 10  \\ \cline{2-11} 
                          & \multirow{4}{*}{2000} & 1000  & 162.35   & 9.24    & 0.12 & 0.02   & 17.57 & 5.98 & 0   & 100 \\
                          &                       & 2000  & 320.72   & 10.89   & 0.26 & 0.08   & 29.46 & 3.35 & 0   & 10  \\
                          &                       & 5000  & 811.04   & 44.56   & 0.58 & 0.18   & 18.20 & 3.19 & 0   & 10  \\
                          &                       & 10000 & 1625.59  & 93.57   & 1.16 & 0.36   & 17.37 & 3.22 & 0   & 10  \\ \cline{2-11} 
                          & \multirow{4}{*}{Mini} & 1000  & 189.59   & 9.24    & 0.29 & 0.06   & 20.52 & 4.96 & 0   & 100 \\
                          &                       & 2000  & 363.14   & 10.89   & 0.43 & 0.13   & 33.36 & 3.36 & 0   & 10  \\
                          &                       & 5000  & 946.20   & 44.56   & 1.12 & 0.46   & 21.23 & 2.41 & 0   & 10  \\
                          &                       & 10000 & 1896.21  & 93.57   & 2.43 & 1.00   & 20.27 & 2.44 & 0   & 10  \\ \hline
\multirow{12}{*}{QMKP}    & \multirow{4}{*}{1000} & 1000  & 1450.09  & 1283.53 & 0.07 & 1.20   & 1.13  & 0.06 & 100 & 100 \\
                          &                       & 2000  & 0.00     & 1967.41 & 0.16 & 3.74   & 0.00  & 0.04 & 10  & 0   \\
                          &                       & 5000  & 14.19    & 0.00    & 0.42 & 10.01  & inf   & 0.04 & 10  & 0   \\
                          &                       & 10000 & 0.00     & 0.00    & 1.13 & 10.01  & /     & 0.11 & 10  & 0   \\ \cline{2-11} 
                          & \multirow{4}{*}{2000} & 1000  & 1159.15  & 1283.53 & 0.09 & 1.22   & 0.90  & 0.07 & 100 & 0   \\
                          &                       & 2000  & 2285.17  & 1967.41 & 0.16 & 3.70   & 1.16  & 0.04 & 10  & 10  \\
                          &                       & 5000  & 8976.34  & 0.00    & 0.58 & 10.01  & inf   & 0.06 & 10  & 10  \\
                          &                       & 10000 & 0.00     & 0.00    & 1.10 & 10.01  & /     & 0.11 & 10  & 0   \\ \cline{2-11} 
                          & \multirow{4}{*}{Mini} & 1000  & 1465.00  & 1283.53 & 0.07 & 1.16   & 1.14  & 0.06 & 100 & 100 \\
                          &                       & 2000  & 2291.27  & 1967.41 & 0.15 & 3.73   & 1.16  & 0.04 & 10  & 10  \\
                          &                       & 5000  & 11418.56 & 0.00    & 0.44 & 10.01  & inf   & 0.04 & 10  & 10  \\
                          &                       & 10000 & 0.00     & 0.00    & 1.07 & 10.01  & /     & 0.02 & 10  & 0   \\ \hline
\end{tabular}%
}
\end{table}

For RandQCP, our method achieves significantly better objective values than Gurobi for each problem instance. Although our runtime is several times longer than Gurobi's, both methods can obtain initial feasible solutions quickly. Furthermore, our models that are trained on Mini and 1000-scale datasets perform well on larger-scale test problems without significant degradation in objective values, even outperforming Gurobi. This demonstrates the strong generalization capability of our framework and validates the feasibility of training our model on small-scale problems and applying it to larger-scale problems. Additionally, models trained on Mini-scale problems perform best, likely due to the larger training dataset, suggesting that increasing the dataset size could further improve model performance.

For QMKP, the model trained on 1000-scale problems outperform Gurobi on 1000-scale test problems, achieving better initial objective values in just 6\% of the time Gurobi took. However, on larger-scale problems, both our method and Gurobi produce initial objective values of 0, likely because heuristic algorithms can quickly find trivial feasible solutions. Since 0 as an objective value is meaningless, we do not make further comparisons. The model trained on Mini-scale problems outperforms Gurobi on 1000-scale and 2000-scale problems, achieving better objective values in just 4\% and 5\% of the time Gurobi takes, respectively. Moreover, on 5000-scale problems, Gurobi produces an objective value of 0 within 10 seconds, whereas our model achieves the best objective value in just 0.44 seconds.

\subsubsection{Parallel Neighborhood Optimization}
\label{appendix:parallel_neighborhood_optimization}

We compared our method with Gurobi and SCIP in the Parallel Neighborhood Optimization stage by running the three methods for a fixed time given the same initial solution. We ran our method on each instance 5 times; we ran Gurobi and SCIP on each instance 1 time. Since our method achieved close objective values among different runs on a single instance, we take the mean as the objective value. We then report the mean and std of the objective value on 10 instances of each method in Table \ref{tab:ablation_lns}. Please refer to supplementary materials for complete raw data.

\begin{table}[H]
    \centering
    \small
    \renewcommand{\arraystretch}{1.1}
    \caption{Comparison of objective values for Parallel Neighborhood Optimization. ``↑" indicates that our method outperforms both Gurobi and SCIP. \textbf{Bold} entries represent the best results.}
    \label{tab:ablation_lns}
    \begin{tabular}{cccccccccc}
        \toprule
        \multirow{2}{*}{Size} & \multirow{2}{*}{Method} & \multicolumn{2}{c}{RandQCP} & & \multicolumn{2}{c}{QMKP}\\\cmidrule{3-4}\cmidrule{6-7}
        & & mean & std & & mean & std\\
        \midrule
        \multirow[c]{4}{*}{5000} & Ours-50\%G & 1495.38↑ & 10.64 & & \textbf{13793.90↑} & 50.01 \\
         & Ours-50\%S & \textbf{1497.16↑} & 10.42 & & 13690.29↑ & 54.62 \\
         & Gurobi & 1428.46 & 10.36 & & 10301.77 & 357.51 \\
         & SCIP & 777.26 & 226.67 & & 8880.72 & 70.70 \\
        \midrule
        \multirow[c]{4}{*}{10000} & Ours-50\%G & 2982.74↑ & 24.42 & & \textbf{14049.71↑} & 45.13 \\
         & Ours-50\%S & \textbf{2987.16↑} & 23.48 & & 14024.71↑ & 55.77 \\
         & Gurobi & 2836.50 & 28.96 & & 13940.55 & 59.86 \\
         & SCIP & 1408.09 & 18.68 & & 10198.27 & 89.47 \\
        \bottomrule
    \end{tabular}    
\end{table}

For the RandQCP problem, our method consistently outperforms both Gurobi and SCIP across all scales. Specifically, the objective values obtained by our method are significantly higher than those obtained by the other solvers. For the 5000-scale problem, our model (Ours-50\%S) achieves an objective value of 1497.16, compared to Gurobi's 1428.46 and SCIP's 777.26. Similarly, for the 10000-scale problem, our model (Ours-50\%S) achieves an objective value of 2987.16, while Gurobi achieves 2836.50 and SCIP only 1408.09. In terms of standard deviation, our method shows consistent results with low variability, indicating the robustness of our approach. 

For the QMKP problem, our method again outperforms the other solvers, especially on larger scales. In the 5000-scale problem, our model (Ours-50\%G) achieves an objective value of 13793.90, while Gurobi and SCIP achieve 10301.77 and 8880.72, respectively. In the 10000-scale problem, our model (Ours-50\%S) achieves an objective value of 14024.71, compared to Gurobi's 13940.55 and SCIP's 10198.27. Moreover, our method's low standard deviations in objective values for QMKP problems indicate the stability and reliability of our framework. For the 5000-scale problem, our model (Ours-50\%S) shows a standard deviation of 54.62, which is lower than that of Gurobi (357.51) and SCIP (70.70).

In summary, our Parallel Neighborhood Optimization method demonstrates superior performance in achieving better objective values with low variability across different problem scales. This validates the efficiency and robustness of our approach, making it a promising solution for large-scale QCQP problems.

\subsection{Experiments on QPLIB}
\label{appendix:QPLIB}

QPLIB is a combination of various mixed problems and datasets, as curated by \cite{furiniQPLIBLibraryQuadratic2019}. The authors filtered out similar problems based on criteria such as variable types, constraint types, and the proportion of non-zero elements. As a result, the 453 problems in QPLIB cover a wide range of possible scenarios, including 133 binary problems.

\subsubsection{Training Dataset Construction}

Given that the problems in QPLIB are not ideally suited for training neural networks, we adopted a novel approach. The problems in QPLIB have fewer than 1000 variables, and their constraints and objective functions can be linear or quadratic. Therefore, we formed a mixed training dataset comprising 1000 instances each of 1000-scale RandQCP and QMKP problems. Subsequently, we used a neural network, trained on this mixed dataset, to predict optimal solutions on QPLIB problems, followed by neighborhood search optimization on the initial feasible solutions. The method of generating this mixed dataset is consistent with Appendix \ref{appendix:dataset_generation}.

\subsubsection{Testing Dataset Construction}
Targeting large-scale, time-consuming problems, we selected certain problems from QPLIB. Initially, we solved these problems using unrestricted Gurobi with a gap limit of 10\% and a time limit of 600s. We filtered out problems that could be optimally solved within 100s and those whose objective values remained the same between 100s and 600s. This process resulted in 16 problems, with variable counts ranging from 144 to 676. We then utilized Gurobi as a small-scale solver at 30\% and 50\% scales, with \verb|TimeLimit| set to 100s and \verb|Threads| set to 4, to compare the objective values achieved at 100s.

\subsubsection{Results}

The experimental results are listed in Table \ref{table:qplib}. Among the 16 problems, our method outperforms Gurobi in 14 cases. In terms of the 30\% and 50\% solver scales, our method matches or exceeds Gurobi in 12 and 10 problems, respectively. Notably, the problems in QPLIB are relatively small-scale; our method's advantages are expected to be more pronounced with even larger-scale problems.

\begin{table}[H]
    \small
    \caption{Comparison of objective values with Gurobi within the same running time on QPLIB. Ours-30\%G and Ours-50\%G mean the scale-limited versions of Gurobi which limit the variable proportion $\alpha$ to 30\% and 50\% respectively. “↑” means the result is equal to or better than the baseline.“n\_var” is the number of variables. “No.” is the QPLIB number. “Sense” is the objective sense.}
    \label{table:qplib}
    \centering
    \begin{tabular}{@{}cccrrr@{}}
        \toprule
        \multirow{2}{*}{n\_var} & \multirow{2}{*}{No.} & \multirow{2}{*}{Sense} & \multicolumn{3}{c}{Objective value} \\ 
        \cmidrule(l){4-6} 
        & & & \multicolumn{1}{c}{Ours-30\%G} & \multicolumn{1}{c}{Ours-50\%G} & \multicolumn{1}{c}{Gurobi} \\ 
        \midrule
        144 & 3402 & min & 230704.0↑           & \textbf{224416.0↑}  & 230704.0          \\
        150 & 5962 & max & \textbf{6962.0↑}    & 6343.0↑             & 5786.0            \\ 
        182 & 3883 & min & \textbf{-788.0↑}    & \textbf{-788.0↑}    & -782.0            \\
        190 & 2067 & min & \textbf{3311060.0↑} & 3441020.0           & 3382980.0         \\
        252 & 2017 & min & -21544.0            & \textbf{-22984.0↑}  & -22584.0          \\
        253 & 2085 & min & 8154640.0           & \textbf{7717850.0↑} & 7885860.0         \\
        275 & 2022 & min & -22000.5↑           & \textbf{-22716.0↑}  & -21514.5          \\
        300 & 3841 & min & -1628.0↑            & \textbf{-1690.0↑}   & -1594.0           \\
        324 & 2036 & min & -28960.0↑           & \textbf{-30480.0↑}  & -28260.0          \\
        324 & 2733 & min & \textbf{5358.0↑}    & \textbf{5358.0↑}    & 5376.0            \\
        435 & 3860 & min & -13820.0            & -16331.0            & \textbf{-16590.0} \\
        462 & 3752 & min & -1075.0             & -1279.0             & \textbf{-1299.0}  \\
        484 & 2957 & min & 3616.0↑             & \textbf{3604.0↑}    & 3788.0            \\
        528 & 3584 & min & \textbf{-13090.0↑}  & -18989.0            & -15525.0          \\
        595 & 2315 & min & \textbf{-13552.0↑}  & -22680.0            & -17952.0          \\
        676 & 3347 & min & \textbf{3825111.0↑} & 3828653.0           & 3826800.0         \\
        \bottomrule
    \end{tabular}%
\end{table}

\section{Algorithms in Parallel Neighborhood Optimization}
\label{appendix:algorithms}

Parallel Neighborhood Optimization is a crucial component of our framework, enabling it to exhibit strong convergence performance and quickly obtain high-quality feasible solutions within a short time. By applying our Q-Repair algorithm, commonly used neighborhood optimization methods for linear problems can be extended to QCQPs. Parallel Neighborhood Optimization in this paper comprises neighborhood partitioning, neighborhood search, and neighborhood crossover. These processes are iterated sequentially to optimize the objective function value repeatedly until the time limit is reached. This section provides a detailed description of the Parallel Neighborhood Optimization process.

\subsection{The Q-Repair Algorithm}
\label{appendix:Q-Repair}

When performing fixed-radius neighborhood optimization for QCQPs, there may be some constraints that are certainly infeasible. The Q-Repair algorithm, based on the McCormick relaxation, quickly identifies these infeasible constraints and reintroduces the correlative variables into the neighborhood, thus rendering the constraints feasible once again.

The first step is using McCormick relaxation (Section \ref{sec:McCormick_relaxation}) to transform quadratic terms in the constraints into linear terms. Then the constraints in QCQPs can be transformed into linear constraints as shown in Equation \ref{eq:qrepair}. Since the repair algorithm is independent of the objective function and only involves constraints, the common Repair algorithm \citep{ye2023gnn} for linear problems can be used to solve this problem after that. Assume that $Q^i_{jk}\geq 0$ and the constraints are all less-equal since the case with negative coefficients and more-equal can be easily generalized. The Q-Repair algorithm is shown in Algorithm \ref{alg:qrepair}.

\begin{algorithm}[H]
    \small
    \caption{Q-REPAIR Algorithm}
    \label{alg:qrepair}
    \begin{algorithmic}[1]
        \REQUIRE The set of fixed variables $\mathcal{F}$, the set of unfixed variables $\mathcal{U}$, the current solution $\vx$, the coefficients of the given QCQP $\{\mQ,\vr,\vb,\vl,\vu\}$\\
        $n\leftarrow$ the number of decision variables\\
        $m\leftarrow$ the number of constraints
        \ENSURE $\mathcal{F},\mathcal{U}$
        \FOR{$i\leftarrow 1$ \TO $n$}
        \IF{$\text{the i-th variable}\in \mathcal{F}$}
            \STATE $cur\_u_i \leftarrow x_i$
            \STATE $cur\_l_i \leftarrow x_i$
        \ELSE
            \STATE $cur\_u_i \leftarrow u_i$
            \STATE $cur\_l_i \leftarrow l_i$
        \ENDIF
        \ENDFOR
        \FOR{$i, j\leftarrow 1$ \TO $n$}
        \STATE $u_{ij},l_{ij} \leftarrow$ \verb+MCCORMICK_RELAXATION+$(u_i,u_j,l_i,l_j)$
        \STATE $cur\_u_{ij},cur\_l_{ij} \leftarrow$ \verb+MCCORMICK_RELAXATION+$(cur\_u_i,cur\_u_j,cur\_l_i,cur\_l_j)$
        \ENDFOR
        \FOR{$i\leftarrow 1$ \TO $m$}
        \STATE $N \leftarrow 0$
        \FOR{$j\leftarrow 1$ \TO $n$}
            \STATE $N\leftarrow N + r^i_jcur\_l_j$
        \ENDFOR
        \FOR{$j, k\leftarrow 1$ \TO $n$}
            \STATE $N\leftarrow N+Q^i_{jk}cur\_l_{jk}$
        \ENDFOR
        \IF{$N>b_i$}
            \FOR{term in constraint $i$}
                \IF{is a linear term with the j-th variable $\in \mathcal{F}$}
                    \STATE remove the j-th variable from $\mathcal{F}$
                    \STATE add the j-th variable into $\mathcal{U}$
                    \STATE $N\leftarrow N - r^i_jcur\_l_j$
                    \STATE $N\leftarrow N + r^i_jl_j$
                \ELSIF{is a quadratic term with the j-th and k-th variable, and at least one of them $\in \mathcal{F}$}
                    \STATE remove them from $\mathcal{F}$
                    \STATE add them into $\mathcal{U}$
                    \STATE $N\leftarrow N - Q^i_{jk}cur\_l_{jk}$
                    \STATE $N\leftarrow N + Q^i_{jk}l_{jk}$
                \ENDIF
                \IF{$N \leq b_i$}
                    \STATE BREAK
                \ENDIF
            \ENDFOR
        \ENDIF
        \ENDFOR
    \end{algorithmic}
 \end{algorithm}

\subsection{Neighboorhood Partition}
\label{appendix:neighborhood_partition}

Neighborhood search refers to fixing a subset of decision variables to their current values and then exploring the search space for the remaining variables. The key to the effectiveness of neighborhood search is the quality of neighborhood partitioning. The neighborhood partition in the proposed framework is an adaptive partitioning based on the density of variables and the size of the neighborhood. It distinguishes the density of the problem by identifying the number of variables contained in each constraint. The number of neighborhoods is determined by the number of variables, neighborhood partition strategy, and neighborhood size.

For dense problems (i.e., one constraint with many variables), the proposed framework employs a variable-based random neighborhood partition strategy. The decision variables are randomly shuffled and added to the neighborhood until the upper limit of the neighborhood size is reached. The number of neighborhoods is equal to the number of variables divided by the neighborhood size. This neighborhood partition strategy ensures that each variable appears in only one neighborhood. Even in dense problems, it guarantees that the solution speed is not compromised due to too many neighborhoods.

For sparse problems (i.e., one constraint with few variables), the proposed framework uses an ACP-based neighborhood partition strategy (\cite{ye2023adaptive}). The constraints are randomly shuffled, and the variables within the constraints are sequentially added to the neighborhood until the upper limit of the neighborhood size is reached. This neighborhood partition strategy increases the likelihood of variables within the same constraint being in the same neighborhood, thereby reducing the probability of the constraint being bound to be violated. More details are shown in Algorithm \ref{alg:neighborhood_partition}.

\begin{algorithm}[H]
    \small
    \caption{ACP-Based Neighborhood Partition}
    \label{alg:neighborhood_partition}
    \begin{algorithmic}[1]
        \REQUIRE A QCQP instance, the number of variables $n$, the number of constraints $m$, the max size of the neighborhood $s_{\text{max}}$
        \ENSURE A set of neighborhoods $\mathcal{N} = \{\mN_1,\mN_2,\dots\}$, the number of neighborhoods $l$
        \STATE Randomly shuffle the order of constraints.
        \STATE $\text{num\_all} \leftarrow 0$
        \FOR{$i\leftarrow 1$ \TO $m$}
        \STATE Initialize $var\_list_i$
        \FOR{$j\leftarrow 1$ \TO $n$}
            \IF{The i-th decision variable is in constraint $i$}
                \STATE add $j$ into $var\_list_i$
                \STATE $\text{num\_all} \leftarrow \text{num\_all} + 1$
            \ENDIF
        \ENDFOR
        \ENDFOR
        \STATE $l \leftarrow \text{num\_all} \, \backslash \, s_{\text{max}}$
        \STATE Initialize $l$ neighborhoods
        \STATE $\text{ID} \leftarrow 1$
        \STATE $s \leftarrow 0$
        \FOR{$i\leftarrow 1$ \TO $m$}
        \FOR{j $\in var\_list_i$}
            \IF{$s < s_{\text{max}}$}
                \STATE add the i-th decision variable into $\mN_{\text{ID}}$
                \STATE $s \leftarrow s+1$
            \ELSE
                \STATE add $\mN_{\text{ID}}$ into $\mathcal{N}$
                \STATE $\text{ID} \leftarrow \text{ID} + 1$
                \STATE $s \leftarrow 1$
                \STATE add the i-th decision variable into $\mN_{\text{ID}}$
            \ENDIF
        \ENDFOR
        \ENDFOR
    \end{algorithmic}
\end{algorithm}

\subsection{Neighboorhood Search and Crossover}
\label{appendix:neighborhood_search_and_crossover}

To prevent the neighborhood search from getting trapped in local optimal, neighborhood crossover is required. The neighborhood crossover used in this paper takes place between two neighborhoods, and the details are provided in Algorithm \ref{alg:crossover}. \verb+Q-REPAIR+ constructs a new search neighborhood to prevent the problem from becoming infeasible after neighborhood crossover. \verb+SEARCH+ means using a small-scale solver to search for a better solution in a specific neighborhood.

Integrating the above steps, the complete Parallel Neighborhood Optimization process is shown in Algorithm \ref{alg:neighborhood_optimization}. The first step is performing neighborhood partitioning and then solving all neighborhoods in parallel using small-scale solvers. After crossing over the neighborhoods pairwise, the best solution is selected as the initial solution for the next round of neighborhood optimization. This process is repeated until the time limit is reached.

\begin{algorithm}[H]
    \small
    \caption{Neighboorhood Crossover}
    \label{alg:crossover}
    \begin{algorithmic}[1]
        \REQUIRE The decision variable set $\mathcal{X}$, neighborhoods $\mN_1,\mN_2$, neighborhood search solution $\vx^1,\vx^2$, the number of variables $n$
        \ENSURE Crossover solution $\vx'$
        \STATE // Assume that $\vx^1$ is better than $\vx^2$
        \FOR{$i\leftarrow 1$ \TO $n$}
        \IF{The i-th decision variable $\in \mN_1$}
            \STATE $\vx'_i \leftarrow \vx^1_i$
        \ELSE
            \STATE $\vx'_i \leftarrow \vx^2_i$
        \ENDIF
        \ENDFOR
        \STATE $\mathcal{F,U}\leftarrow$ \verb+Q-REPAIR+($\mathcal{X},\emptyset,\vx'$)
        \STATE $\mathcal{X}'\leftarrow$ \verb+SEARCH+($\mathcal{F,U},\vx'$)
    \end{algorithmic}
\end{algorithm}

\begin{algorithm}[H]
    \small
    \caption{Parallel Neighborhood Optimization}
    \label{alg:neighborhood_optimization}
    \begin{algorithmic}[1]
        \REQUIRE Initial feasible solution $\vx$, the number of variables $n$, the max variable proportion $\alpha_{\text{ub}}$, time limit $t$
        \ENSURE Optimization solution $\vx$
        \STATE $\mathcal{N},l\leftarrow$ \verb+PARTITION+($s_{\text{max}} = \alpha_{\text{ub}}\cdot n$)
        \STATE // Do the next step in parallel
        \FOR{$\mathcal{F}_i,\mathcal{U}_i \in \mathcal{N}$ }
            \STATE $\vx_i \leftarrow$ \verb+SEARCH+($\mathcal{F}_i, \mathcal{U}_i, \vx$)
        \ENDFOR
        \STATE // Do the next step in parallel
        \FOR{$i\leftarrow 1$ \TO $\floor{l/2}$}
            \STATE $\vx'_i \leftarrow$ \verb+CROSSOVER+($\vx_{2i-1},\vx_{2i},\mN_{2i-1},\mN_{2i}$)
        \ENDFOR
        \STATE $\vx\leftarrow$ the best solution among $\vx'_i$
        \IF{time limit reaches}
            \RETURN $\vx$
        \ELSE
            \STATE Restart from row 1 with $\vx$
        \ENDIF
    \end{algorithmic}
\end{algorithm}

\section{Limitations, Broader Impacts and Future Directions}
\label{appendix:limitation_impacts_future_directions}

In this section, we provide a discussion of limitations, broader impacts, and future directions of our framework.

\subsection{Discussion of Limitations}
\label{appendix:discussion_limitations}

Limitations of our paper include:

\begin{enumerate}[left=0pt, topsep=0pt, itemsep=0pt]
    \item Our framework is designed for large-scale QCQP problems. However, due to the lack of publicly available datasets and the limited number of problems in QPLIB, we had to generate large-scale QCQP problems to test the effectiveness of our framework.
    \item We have proven the equivalence of UniEGNN and our hypergraph representation to the Interior-Point Method (IPM) only for QPs. This equivalence has not been extended to more general nonlinear cases. In the appendix, we discuss the nonlinear interior-point method to facilitate future work in extending our framework to general nonlinear problems.
    \item While our framework is suitable for general QCQP problems, due to the lack of standard problem formulations and datasets for general integer variables, we have only tested our framework on problems with binary variables.
\end{enumerate}

\subsection{Discussion of Broader Impacts}
\label{appendix:discussion_broader_impacts}

The optimization model introduced in this paper primarily aims to solve large-scale QCQPs efficiently. This has several broader impacts:

\begin{itemize}[left=0pt, topsep=0pt, itemsep=0pt]
    \item \textbf{Positive Impacts:}
        \begin{itemize}[left=0pt, topsep=0pt, itemsep=0pt]
            \item \textbf{Enhanced Optimization Solutions:} The proposed model can improve the efficiency and accuracy of solving complex optimization problems, benefiting various fields such as logistics, finance, engineering, and operations research.
            \item \textbf{Research and Development:} The model can serve as a foundation for further academic research, leading to advancements in optimization techniques and their applications.
            \item \textbf{Industrial Applications:} Improved optimization methods can enhance decision-making processes in industries, leading to cost savings, resource optimization, and better overall performance.
        \end{itemize}
    \item \textbf{Negative Impacts:}
        \begin{itemize}[left=0pt, topsep=0pt, itemsep=0pt]
            \item \textbf{Misuse in Competitive Scenarios:} While unlikely, there is a potential risk that the model could be misused in highly competitive industries to gain unfair advantages or manipulate market conditions.
            \item \textbf{Over-Reliance on Automated Systems:} Increased reliance on automated optimization solutions may lead to reduced human oversight, potentially causing issues if the model produces incorrect or biased results.
        \end{itemize}
\end{itemize}

Overall, the societal impacts of this work are primarily positive, with the potential to advance the field of optimization and its applications significantly. The negative impacts are limited and can be mitigated with proper oversight and ethical guidelines for the use of such models.

\subsection{Future Directions}
\label{appendix:future_directions}

Based on the limitations identified, several future directions can be pursued:

\begin{enumerate}[left=0pt, topsep=0pt, itemsep=0pt]
    \item \textbf{Expanding Dataset Availability:} Efforts should be made to develop and release publicly available datasets for large-scale QCQP problems. This would facilitate more comprehensive testing and validation of optimization frameworks like NeuralQP.
    \item \textbf{Extending to Nonlinear Optimization:} Future research should aim to extend the equivalence of UniEGNN and IPM to more general nonlinear optimization problems. Building on the discussion of nonlinear IPMs provided in the appendix, this extension could enhance the applicability of NeuralQP to a broader range of optimization challenges.
    \item \textbf{Testing on General Integer Variables:} Developing standard problem formulations and datasets for problems with general integer variables will enable thorough evaluation and potential improvements to the framework's versatility.
\end{enumerate}

\end{document}